\documentclass[11pt,twoside,final]{scrartcl}
\usepackage[english]{babel}
\usepackage{url,a4wide}
\usepackage{amsmath, amsfonts, amssymb, amsthm, mathtools, nicefrac}
\usepackage{todonotes}
\usepackage{enumitem}
\usepackage{graphicx}
\graphicspath{{pics/}}
\usepackage{xcolor}
\usepackage[skip=1pt, font=normalsize]{subcaption}
\usepackage{hyperref} % provides \url command for bibtex and links to jump within documents
\hypersetup{plainpages=false, colorlinks, linkcolor=black, citecolor=black, urlcolor=blue,
pdftitle={Some remarks on regularized Shannon sampling formulas},
pdfauthor={Melanie Kircheis, Daniel Potts, Manfred Tasche},
pdfstartview={FitBH}}
\usepackage{placeins}
\usepackage{xurl}
\allowdisplaybreaks

\newcommand{\e}{\mathrm e}
\renewcommand{\i}{\mathrm i}

\newcommand{\sinc}{\mathrm{sinc}}

\newcommand{\R}{\mathbb R}
\newcommand{\C}{\mathbb C}
\newcommand{\Z}{\mathbb Z}
\newcommand{\N}{\mathbb N}

\newcommand{\new}[1]{\textcolor{black}{ #1}}

\newtheorem{theorem}{Theorem}[section]
\newtheorem{corollary}[theorem]{Corollary}
\newtheorem{lemma}[theorem]{Lemma}
\newtheorem{alg}[theorem]{Algorithm}
\newtheorem{definition}[theorem]{Definition}
\newtheorem{example}[theorem]{Example}
\newtheorem{remark}[theorem]{Remark}
\newtheorem{conjecture}[theorem]{Conjecture}

\newcommand{\bend}{\hspace*{0ex} \hfill \hbox{\vrule height
	1.5ex\vbox{\hrule width 1.4ex \vskip 1.4ex\hrule width 1.4ex}\vrule
	height 1.5ex}}
\renewcommand{\qedsymbol}{\rule{1.5ex}{1.5ex}}

\newenvironment{Lemma}{\goodbreak\begin{lemma}}{\end{lemma}}
\newenvironment{Theorem}{\goodbreak\begin{theorem}}{\end{theorem}}
\newenvironment{Remark}{\goodbreak\begin{remark}\upshape}{\bend\end{remark}}

\newenvironment{Example}{\goodbreak\begin{example}\upshape}{\bend\end{example}}
\newenvironment{Conjecture}{\goodbreak\begin{conjecture}\upshape}{\end{conjecture}}
\newenvironment{algorithm}[1]{\goodbreak~\begin{alg}[#1]~\vspace{-9pt}~\\
		\rule{\linewidth}{0.5pt}~\\}{\vspace{-9pt}~\\
		\rule{\linewidth}{0.5pt}~\end{alg}}

\numberwithin{equation}{section}
\numberwithin{table}{section}
\numberwithin{figure}{section}

\renewcommand{\mathbf}[1]{\ensuremath{\boldsymbol{#1}}}

\title{Some remarks on regularized Shannon sampling formulas}
\author{
	Melanie Kircheis\footnotemark[1] \and
	Daniel Potts\footnotemark[4] \and
	Manfred Tasche\footnotemark[3]
}

\date{}
\begin{document}
\maketitle

\begin{abstract}

The fast reconstruction of a bandlimited function from its sample data is an essential problem in signal processing.
In this paper, we consider the widely used Gaussian regularized Shannon sampling formula in comparison to regularized Shannon sampling formulas employing alternative window functions, \new{such as} the $\sinh$-type window function and the continuous Kaiser--Bessel window function.
It is shown that the approximation errors of these regularized Shannon sampling formulas possess an exponential decay with respect to the truncation parameter.
\new{The main focus of this work is to address minor gaps in the preceding papers~\cite{KiPoTa22,KiPoTa24} and rigorously prove assumptions that were previously based solely on numerical tests.}
In doing so, we demonstrate that the \new{$\sinh$-type regularized Shannon sampling formula has the same exponential decay as the continuous Kaiser--Bessel regularized Shannon sampling formula, but both have twice the exponential decay of the Gaussian regularized Shannon sampling formula.}
Additionally, numerical experiments illustrate the theoretical results.
\medskip
	
\smallskip
	
AMS \emph{Subject Classifications}:
62D05, 42A10, 41A25, 94A20.
\end{abstract}

\footnotetext[1]{Corresponding author: melanie.kircheis@math.tu-chemnitz.de, Chemnitz University of Technology, Faculty of Mathematics, D--09107 Chemnitz, Germany}
\footnotetext[4]{potts@mathematik.tu-chemnitz.de, Chemnitz University of Technology, Faculty of Mathematics, D--09107 Chemnitz, Germany}
\footnotetext[3]{manfred.tasche@uni-rostock.de, University of Rostock, Institute of Mathematics, D--18051 Rostock, Germany}

\section{Introduction}

In signal processing, the fast reconstruction of a bandlimited function from its sample data is of fundamental importance. A function~\mbox{$f \in L^2(\mathbb R) \cap C(\mathbb R)$} is called
\emph{bandlimited} with \emph{bandwidth}~\mbox{$\delta > 0$}, if its Fourier transform
\begin{align}
\label{eq:fourier_trafo}
({\mathcal F}f)(\omega) = {\hat f}(\omega)\coloneqq \frac{1}{\sqrt{2 \pi}}\, \int_{\mathbb R} f(t)\,{\mathrm e}^{-{\mathrm i}t \omega}\,{\mathrm d}t\,, \quad \omega \in \mathbb R\,,
\end{align}
vanishes for all~\mbox{$|\omega| \geq \delta$}. %By the oversampling condition~\mbox{$\delta < \pi$},
For such a bandlimited function with~\mbox{$\delta \in (0,\,\pi]$} the famous Shannon sampling theorem, see~\cite{Whittaker, Kotelnikov, Shannon49}, states that
\begin{align}
\label{eq:Shannonseries}
f(t) = \sum_{k\in \mathbb Z} f(k)\,\sinc (t-k)\,, \quad t \in \mathbb R\,,
\end{align}
where
\begin{align}
\label{eq:sinc}
\sinc (t) \coloneqq \begin{cases} \frac{\sin(\pi t)}{\pi t} &\colon t \in \mathbb R \setminus \{0\}\,,\\
1 &\colon t=0\,,
\end{cases}
\end{align}
denotes the \emph{cardinal sine function}. It is known that the Shannon sampling series~\eqref{eq:Shannonseries} converges absolutely and uniformly on whole~$\mathbb R$.
However, the practical use of~\eqref{eq:Shannonseries} is limited, since its evaluation requires infinitely many samples %of~$f$
and its truncated version is not a good approximation due to the slow decay of the cardinal sine function, see~\cite{Ja66}.
In addition to this rather poor convergence, it is known, see~\cite{Fe92a,Fe92b,DDeV03}, that in the presence of noise in the samples~\mbox{$f(k)$}, \mbox{$k\in\Z$}, of a bandlimited function~\mbox{$f \in L^2(\mathbb R) \cap C(\mathbb R)$} the convergence of Shannon sampling series~\eqref{eq:Shannonseries} may even break	down completely.
Therefore, it was proposed to consider the regularization of the Shannon sampling series with a suitable window function.
Note that many authors such as~\cite{D92, Nat86, Rap96, Par97, StTa06} used window functions in the frequency domain, but the recent study~\cite{KiPoTa24} has shown that it is much more beneficial to employ a window function in the spatial domain, cf.~\cite{Q03, Q04, StTa06, MXZ09, LZ16, CZ19, KiPoTa22}.
In the following, a \emph{window function}~\mbox{$\varphi:\,\R \to [0,\,1]$} is an even function in~\mbox{$L^2(\mathbb R) \cap C(\mathbb R)$} which decreases on~\mbox{$[0,\,\infty)$} and fulfills~\mbox{$\varphi(0) = 1$}.
By~\mbox{${\mathbf 1}_{[-m,\,m]}$} we denote the \emph{characteristic function} of the interval~\mbox{$[-m,\,m]$} with~\mbox{$m \in \N \setminus \{1\}$}, i.\,e., the function
\begin{align*}
{\mathbf 1}_{[-m,\,m]}(t) \coloneqq \begin{cases} 1 & \colon t \in [-m,\,m]\,, \\
0 & \colon t \in \R \setminus [-m,\,m]\,.
\end{cases}
\end{align*}
In this paper, we assume that the bandwidth~$\delta$ of~$f$ fulfills the so-called \emph{oversampling condition}~\mbox{$0<\delta<\pi$}.
Then we recover~$f$ by the \emph{\mbox{$\varphi$-reg}ularized Shannon sampling formula}
\begin{align}
\label{eq:regShannonformula}
\big(R_{\varphi,m}f\big)(t) \coloneqq \sum_{k\in \Z} f(k)\, \sinc(t-k)\,\varphi(t-k)\,{\mathbf 1}_{[-m,\,m]}(t-k)\,, \quad t\in \R\,,
\end{align}
where~\mbox{$m \in \N \setminus \{1\}$} is the so-called \emph{truncation parameter}. In doing so, we consider the following window functions~\mbox{$\varphi: \, \R \to [0,\,1]$}.

\begin{Remark}
\label{Example:windows}
The most popular window function, see e.\,g.~\cite{Q03, QO05, SchSt07, TaSuMu07, LZ16, CZ19}, is the \emph{Gaussian function}
\begin{align}
\label{eq:Gaussfunction}
\varphi_{\mathrm{Gauss}}(t) \coloneqq \e^{-t^2/(2 \sigma^2)}\,, \quad t \in \mathbb R\,,
\end{align}
with \emph{variance}~\mbox{$\sigma^2 > 0$}.
Note that \new{this} window function \new{is} supported on whole~$\R$.

Here we prefer window functions which are compactly supported on the interval~\mbox{$[-m,\,m]$}, as studied in~\cite{KiPoTa22, KiPoTa24}.
The $\sinh$-\emph{type window function} is defined as
\begin{align}
\label{eq:varphisinh}
	\varphi_{\sinh}(t) \coloneqq
	\begin{cases}
		\frac{1}{\sinh \beta}\, \sinh\Big(\beta\,\sqrt{1-\tfrac{t^2}{m^2}}\,\Big) &\colon t \in [-m,\,m] \,, \\
		0 &\colon t \in \mathbb R\setminus [-m,\,m]  \,,
	\end{cases}
\end{align}
with \emph{shape parameter}~\mbox{$\beta > 0$}, see~\cite{PT21a}. Then the corresponding expression~\eqref{eq:regShannonformula} is termed the $\sinh$-\emph{type regularized Shannon sampling formula}.
The \emph{continuous Kaiser--Bessel window function} is defined as
\begin{align}
\label{eq:varphicKB}
	\varphi_{\mathrm{cKB}}(t) \coloneqq
	\begin{cases}
		\frac{1}{I_0(\beta)-1}\,\big(I_0(\beta \sqrt{1- t^2/m^2}) - 1\big) &\colon t \in [-m,\,m] \,, \\
		0 &\colon t \in \mathbb R\setminus [-m,\,m]  \,,
	\end{cases}
\end{align}
with convenient shape parameter~\mbox{$\beta > 0$}, see~\cite{PT21a}. Then the corresponding expression~\eqref{eq:regShannonformula} is called the \emph{continuous Kaiser--Bessel regularized Shannon sampling formula}.
We remark that these two window functions~\eqref{eq:varphisinh} and~\eqref{eq:varphicKB} are well-studied in the context of the nonuniform fast Fourier transform (NFFT), see e.\,g.~\cite[Section~6]{PPST23} and~\cite{BaMaKl18, Ba20}.
\end{Remark}

Due to the definition of the cardinal sine function~\eqref{eq:sinc} we have~\mbox{$\sinc(n-k)=\delta_{n,k}$} and therefore the regularized Shannon sampling formula~\mbox{$R_{\varphi,m}f$} in~\eqref{eq:regShannonformula} has the \emph{interpolation property}
\begin{align}
\label{eq:interpolprop}
\big(R_{\varphi,m}f\big)(n) = f(n)\,, \quad n \in \Z\,.
\end{align}
Moreover, the use of the characteristic function~\mbox{${\mathbf 1}_{[-m,\,m]}$} in~\eqref{eq:regShannonformula} leads to \emph{localized sampling} of~$f$, i.\,e., the
computation of~\mbox{$\big(R_{\varphi,m}f\big)(t)$} for any~\mbox{$t \in \R \setminus \Z$} requires only~\mbox{$2m$} samples~\mbox{$f(k)$}, where~\mbox{$k \in \Z$} fulfills the condition
\mbox{$|k -t| \leq m$}. Especially, for~\mbox{$t \in (0,\,1)$} we obtain the finite sum
\begin{align*}
\big(R_{\varphi,m}f\big)(t) = \sum_{k=1-m}^m f(k)\, \sinc(t-k)\,\varphi(t-k)\,.
\end{align*}
As in many applications, we use \emph{oversampling} of the given bandlimited function~$f$ with bandwidth~\mbox{$\delta < \pi$}, i.\,e., the function~$f$ is sampled on the integer grid~$\Z$.

In this paper, we focus on the \mbox{$\varphi$-reg}ularized Shannon sampling formulas~\eqref{eq:regShannonformula} for the window functions~\mbox{$\varphi$} given in Remark~\ref{Example:windows}. To compare the corresponding approaches, we present estimates \new{of} the uniform approximation error
\begin{align}
\label{eq:err_approx}
\| f - R_{\varphi,m}f \|_{C_0(\R)} \coloneqq \max_{t\in \R} \big| f(t) - \big(R_{\varphi,m}f\big)(t)\big|\,,
\end{align}
where~\mbox{$C_0(\R)$} denotes the Banach space of continuous functions~\mbox{$g\colon\!\R \to \C$} vanishing as~\mbox{$|t| \to \infty$} \new{equipped with the norm~\mbox{$\|f\|_{C_0(\R)} \coloneqq \max_{t\in \R}|f(t)|$}}.
\new{Primarily, this work concentrates on addressing minor gaps in the preceding papers~\cite{KiPoTa22,KiPoTa24} and rigorously proving the corresponding assumptions that were previously based solely on numerical experiments.}

For this purpose, we initially study the uniform approximation error of general \mbox{$\varphi$-reg}ularized Shannon sampling formulas~\eqref{eq:regShannonformula} in Section~\ref{sec:approx_error}.
Afterwards, we specify our findings for the window functions~$\varphi$ introduced in Remark~\ref{Example:windows}.
In particular, Section~\ref{sec:Gauss} deals with the Gaussian window function~\eqref{eq:Gaussfunction}, while Section~\ref{sec:sinh} is concerned with the $\sinh$-type window function~\eqref{eq:varphisinh} and Section~\ref{sec:cKB} with the continuous Kaiser--Bessel window function~\eqref{eq:varphicKB}.

\section{Approximation error of regularized Shannon sampling formulas \label{sec:approx_error}}

Firstly, we estimate the uniform approximation error of the \mbox{$\varphi$-reg}ularized Shannon sampling formula~\eqref{eq:regShannonformula}, analogously to~\cite[Theorem~3.2]{KiPoTa22} and~\cite[Theorem~4.1]{KiPoTa24}.

\begin{Theorem}
\label{Theorem:approxerror}
Assume that~\mbox{$f\in L^2(\R) \cap C(\R)$} is bandlimited with bandwidth~\mbox{$\delta \in (0,\,\pi)$}. Further let~\mbox{$\varphi:\, \R \to [0,\,1]$} be an even function in~\mbox{$L^2(\R) \cap C(\R)$}
which is decreasing on~\mbox{$[0, \,\infty)$} with~\mbox{$\varphi(0) = 1$},
and let~\mbox{$m \in {\mathbb N}\setminus \{1\}$} be given. \\
Then the \mbox{$\varphi$-reg}ularized Shannon sampling formula~\eqref{eq:regShannonformula} satisfies the error estimate
\begin{align*}
	\| f - R_{\varphi,m}f \|_{C_0(\mathbb R)} \le \big( E_1(m) +  E_2(m) \big) \,\|f\|_{L^2(\mathbb R)} \,, \quad m \in \N \setminus \{1\}\,,
\end{align*}
with the error constants
\begin{align}
	E_1(m) &\coloneqq \max_{\omega \in [-\delta,\,\delta]} \bigg| 1 - \frac{1}{\sqrt{2\pi}}\,\int_{\omega - \pi}^{\omega + \pi} {\hat\varphi}(\tau)\,{\mathrm d}\tau \,\bigg|\,, \label{eq:E1} \\
	E_2(m) &\coloneqq \frac{\sqrt{2}}{\pi\, m}\, \sqrt{\varphi^2(m) + \int_m^{\infty} \varphi^2(t) \,\mathrm{d}t}\,. \label{eq:E2}
\end{align}
\end{Theorem}

\emph{Proof}. (i) Initially, we consider only the case~\mbox{$t \in (0,\,1)$}, where we split the approximation error
\begin{align*}
f(t) - \big(R_{\varphi,m}f\big)(t) = e_1(t) + e_{2,0}(t)\,, \quad t\in (0,\,1)\,,
\end{align*}
into the \emph{regularization error}
\begin{align}
\label{eq:e1}
	e_1(t) \coloneqq f(t) - \sum_{k\in \Z} f(k)\, \sinc(t-k)\,\varphi(t-k)\,, \quad t \in \R\,,
\end{align}
and the \emph{truncation error}
\begin{align}
\label{eq:e20}
e_{2,0}(t) &\coloneqq  \sum_{k\in \Z} f(k)\, \sinc(t-k)\,\varphi(t-k) - \big(R_{\varphi,m}f\big)(t) \nonumber\\
&= \sum_{k \in \Z \setminus \{1-m,\ldots,m\}} f(k)\, \sinc(t-k)\,\varphi(t-k)\,, \quad t \in (0,\,1)\,.
\end{align}
(ii) To estimate the regularization error~\eqref{eq:e1}, we start our study by considering the Fourier transform~\eqref{eq:fourier_trafo} of the function~\mbox{$\varphi \,\sinc$}, i.\,e., the term
\begin{align*}
{\mathcal F}(\varphi\,\sinc)(\omega) = \frac{1}{\sqrt{2\pi}}\, \int_{\R} \varphi(t)\,\sinc (t)\, \e^{-\i \omega t}\, \mathrm{d}t\,.
\end{align*}
Using the convolution property of~$\mathcal F$ in~\mbox{$L^2(\R)$} (see~\cite[Theorem 2.26]{PPST23}), we have
\begin{align*}
{\mathcal F}(\varphi\,\sinc)(\omega) = \big({\hat \varphi}\star ({\mathcal F}\sinc)\big)(\omega) %\\
= \frac{1}{\sqrt{2\pi}}\, \int_{\R} \hat \varphi(\omega - \tau)\,(\mathcal F \sinc) (\tau)\, \mathrm{d}\tau\,,
\end{align*}
and hence by
\begin{align*}
(\mathcal F \sinc) (\tau) = \frac{1}{\sqrt{2\pi}}\,{\mathbf 1}_{[-\pi,\pi]}(\tau)
\end{align*}
we obtain
\begin{align*}
{\mathcal F}(\varphi\,\sinc)(\omega) = \frac{1}{2\pi}\,\int_{\omega - \pi}^{\omega+ \pi} \hat \varphi(\tau)\,\, \mathrm{d}\tau\,.
\end{align*}
Consequently, using the shifting property of~$\mathcal F$, the Fourier transform~\eqref{eq:fourier_trafo} of the shifted function~\mbox{$\varphi(t-k)\,\sinc (t-k)$} with~\mbox{$k\in \Z$} reads as
\begin{align*}
\frac{1}{\sqrt{2\pi}}\,\int_{\R} \varphi(t-k)\,\sinc (t-k)\,\e^{-\i\omega t}\,\mathrm{d}t = \e^{-\i \omega k}\,{\mathcal F}(\varphi\,\sinc)(\omega)
= \frac{1}{2\pi}\,\e^{-\i \omega k}\,\int_{\omega - \pi}^{\omega+ \pi} \hat \varphi(\tau)\,\, \mathrm{d}\tau\,.
\end{align*}
Therefore, the Fourier transform of the regularization error~$e_1$ in~\eqref{eq:e1} has the form
\begin{align}
\label{eq:hate1}
{\hat e}_1(\omega) = {\hat f}(\omega) - \bigg(\frac{1}{2\pi}\, \sum_{k\in \Z} f(k) \,\e^{-\i \omega k}\bigg)\,\int_{\omega - \pi}^{\omega+ \pi} \hat \varphi(\tau)\,\, \mathrm{d}\tau\,.
\end{align}
Note that since the set of shifted cardinal sine functions~\mbox{$\mathrm{sinc}(\cdot - k)$} with~\mbox{$k \in \mathbb Z$} forms an orthonormal system in~\mbox{$L^2(\mathbb R)$}, i.\,e.,
\begin{align*}
\int_{\R} \sinc (t-k)\,\sinc (t-\ell)\,\mathrm{d}t = \delta_{k,\ell}\,, \quad k,\, \ell \in \Z\,,
\end{align*}
and the given function~$f$ can be represented by the Shannon sampling series~\eqref{eq:Shannonseries}, we obtain that
\begin{align}
\label{eq:parseval}
\sum_{k\in \mathbb Z} |f(k)|^2 &= \sum_{k\in \mathbb Z} \sum_{\ell\in \mathbb Z} f(k) \,\overline{f(\ell)} \int_{\R} \sinc (t-k)\,\sinc (t-\ell) \,\mathrm{d}t \notag \\
&= \int_{\R} f(t) \,\overline{f(t)} \,\mathrm{d}t = \| f \|_{L^2(\mathbb R)}^2 < \infty \,,
\end{align}
and thus the series
\begin{align*}
\sum_{k\in \mathbb Z} f(k)\, {\mathrm e}^{- {\mathrm i} \omega k}
\end{align*}
converges in~\mbox{$L^2([-\pi,\pi])$}.
Moreover, since~$f$ is bandlimited with bandwidth~\mbox{$\delta \in (0,\,\pi)$}, we have~\mbox{${\hat f}(\omega) = 0$} for all~\mbox{$\omega \in \R \setminus [- \delta,\,\delta]$},
and thereby the restricted function~\mbox{$\hat f\big|_{[-\pi,\,\pi]}$} belongs to~\mbox{$L^2([-\pi,\pi])$}. Hence, this restricted function possesses the~\mbox{$2\pi$}-periodic
Fourier expansion
\begin{align*}
{\hat f}(\omega) = \sum_{k \in \Z} c_k(\hat f)\,\e^{-\i \omega k}\,, \quad \omega \in [-\pi,\, \pi]\,,
\end{align*}
with the Fourier coefficients
\begin{align*}
c_k(\hat f) &= \frac{1}{2\pi}\,\int_{-\pi}^{\pi} {\hat f}(\tau)\,\e^{\i k \tau}\,\mathrm{d}\tau =  \frac{1}{2\pi}\,\int_{\R} {\hat f}(\tau)\,\e^{\i k \tau}\,\mathrm{d}\tau%\\&
= \frac{1}{\sqrt{2\pi}}\,f(k)\,, \quad k \in \Z\,,
\end{align*}
by inverse Fourier transform. In other words, the function~$\hat f$ can be represented in the form
\begin{align}
\label{eq:hatf}
\hat f(\omega) = {\hat f}(\omega)\,{\mathbf 1}_{[-\delta,\delta]}(\omega) = \frac{1}{\sqrt{2\pi}}\,\bigg(\sum_{k\in \Z} f(k)\,\e^{-\i k \omega}\bigg)\,{\mathbf 1}_{[-\delta,\delta]}(\omega)\,, \quad \omega \in \R\,.
\end{align}
Introducing the auxiliary function
\begin{align*}
\Delta_{\varphi}(\omega)\coloneqq {\mathbf 1}_{[-\delta,\delta]}(\omega) - \frac{1}{\sqrt{2\pi}}\,\int_{\omega - \pi}^{\omega+ \pi} \hat \varphi(\tau)\,\, \mathrm{d}\tau\,, \quad \omega \in \R\,,
\end{align*}
we see by inserting~\eqref{eq:hatf} into~\eqref{eq:hate1} that
\begin{align*}
{\hat e}_1(\omega) = {\hat f}(\omega)\,\Delta_{\varphi}(\omega)\,, \quad \omega \in \R\,,
\end{align*}
and thereby
\begin{align*}
\big|{\hat e}_1(\omega)\big| = \big|{\hat f}(\omega)\big|\,\big|\Delta_{\varphi}(\omega)\big|\,, \quad \omega \in \R\,.
\end{align*}
Thus, inverse Fourier transform and the definition~\eqref{eq:E1} yields
\begin{align*}
|e_1(t)| &\leq \frac{1}{\sqrt{2\pi}}\,\int_{\R} \big|{\hat e}_1(\omega)\big|\,\mathrm{d}\omega = \frac{1}{\sqrt{2\pi}}\,\int_{-\delta}^{\delta} \big|{\hat f}(\omega)\big|\,\big|\Delta_{\varphi}(\omega)\big|\,\mathrm{d}\omega\\
&\leq \frac{1}{\sqrt{2\pi}}\, \max_{\omega \in [-\delta,\delta]} \big|\Delta_{\varphi}(\omega)\big|\;\int_{-\delta}^{\delta} \big|{\hat f}(\omega)\big|\,\mathrm{d}\omega%\\
= \frac{1}{\sqrt{2\pi}}\,E_1(m)\,\int_{-\delta}^{\delta} \big|{\hat f}(\omega)\big|\,\mathrm{d}\omega\,.
\end{align*}
By the Cauchy--Schwarz inequality and the Parseval equality~\mbox{$\| \hat f\|_{L^2(\R)} = \| f \|_{L^2(\R)}$} we obtain
\begin{align*}
\int_{-\delta}^{\delta} \big| 1 \cdot {\hat f}(\omega)\big| \, \mathrm{d}\omega &\leq \bigg(\int_{-\delta}^{\delta} 1^2\,\mathrm{d}\omega\bigg)^{1/2}\,\bigg(\int_{-\delta}^{\delta} \big|{\hat f}(\omega)\big|^2\,\mathrm{d}\omega\bigg)^{1/2} %\\ &
= \sqrt{2\delta}\,\| \hat f\|_{L^2(\R)} \leq \sqrt{2\pi}\,\| f \|_{L^2(\R)}\,.
\end{align*}
Consequently, we receive the estimate
\begin{align*}
|e_1(t)| \leq E_1(m)\,\| f \|_{L^2(\R)} \,, \quad t \in \R\,,
\end{align*}
and hence
\begin{align*}
\max_{t \in \R} |e_1(t)| \leq E_1(m)\,\| f \|_{L^2(\R)}\,.
\end{align*}
(iii) Now we estimate the truncation error~\mbox{$e_{2,0}(t)$} for~\mbox{$t\in (0,\,1)$}. By~\eqref{eq:e20} and~\mbox{$\varphi(t) \geq 0$}, we obtain
\begin{align*}
|e_{2,0}(t)| \leq \sum_{k \in \Z \setminus \{1-m,\ldots,m\}} |f(k)|\, |\sinc(t-k)|\,\varphi(t-k)\,, \quad t \in (0,\,1)\,.
\end{align*}
For~\mbox{$t \in (0,\,1)$} and~\mbox{$k \in \Z \setminus \{1-m,\ldots,m\}$}, we estimate
\begin{align*}
 |\sinc(t-k)| \leq \frac{1}{\pi\,|t-k|} \leq \frac{1}{\pi m}\,,
\end{align*}
such that
\begin{align*}
|e_{2,0}(t)| \leq \frac{1}{\pi m}\,\sum_{k \in \Z \setminus \{1-m,\ldots,m\}} |f(k)|\,\varphi(t-k)\,, \quad t \in (0,\,1)\,.
\end{align*}
Then the Cauchy--Schwarz inequality implies
\begin{align*}
|e_{2,0}(t)| \leq \frac{1}{\pi m}\,\bigg(\sum_{k \in \Z \setminus \{1-m,\ldots,m\}} |f(k)|^2\bigg)^{1/2}\,\bigg(\sum_{k \in \Z \setminus \{1-m,\ldots,m\}} \varphi^2(t-k)\bigg)^{1/2}\,, \quad t \in (0,\,1)\,.
\end{align*}
From~\eqref{eq:parseval} it follows that
\begin{align*}
|e_{2,0}(t)| \leq \frac{1}{\pi m}\, \| f\|_{L^2(\R)}\,\bigg(\sum_{k \in \Z \setminus \{1-m,\ldots,m\}} \varphi^2(t-k)\bigg)^{1/2}\,, \quad t \in (0,\,1)\,.
\end{align*}
Since by assumption the window function~$\varphi$ is even and~\mbox{$\varphi\big|_{[0,\infty)}$} decreases, we can estimate the series
\begin{align*}
\sum_{k \in \Z \setminus \{1-m,\ldots,m\}} \varphi^2(t-k) &= \sum_{k=-\infty}^{-m} \varphi^2(t-k) + \sum_{k=m+1}^{\infty} \varphi^2(t-k)\\
&= \sum_{k=m}^{\infty} \varphi^2(t+k) + \sum_{k=m+1}^{\infty} \varphi^2(k-t)\\
&\leq  \sum_{k=m}^{\infty} \varphi^2(k) + \sum_{k=m+1}^{\infty} \varphi^2(k-1) = 2 \,\sum_{k=m}^{\infty} \varphi^2(k) \,, \quad t \in (0,\,1)\,.
\end{align*}
Applying the integral test for convergence of series, we obtain that
\begin{align*}
 2 \,\sum_{k=m}^{\infty} \varphi^2(k) = 2\,\varphi^2(m) + 2\,\sum_{k=m+1}^{\infty} \varphi^2(k) < 2\,\varphi^2(m) + 2\,\int_{m}^{\infty} \varphi^2(t)\,\mathrm{d}t\,.
\end{align*}
Thus, for each~\mbox{$t \in (0,\,1)$} we have by definition~\eqref{eq:E2} that
\begin{align*}
|e_{2,0}(t)| \leq \frac{\sqrt{2}}{\pi m}\, \bigg(\varphi^2(m) + \int_m^{\infty} \varphi^2(t) \,\mathrm{d}t\bigg)^{1/2}\, \| f\|_{L^2(\R)} = E_2(m)\, \| f\|_{L^2(\R)} < \infty\,.
\end{align*}
Furthermore, by the interpolation property~\eqref{eq:interpolprop} of~\mbox{$R_{\varphi,m}f$} we have~\mbox{$e_{2,0}(0) = e_{2,0}(1) = 0$}, such that
\begin{align*}
\max_{t\in [0,1]} |e_{2,0}(t)| \leq E_2(m)\,\| f\|_{L^2(\R)}\,.
\end{align*}
(iv) By the same technique, the error estimate
\begin{align*}
\max_{t \in [n, n+1]} \big| f(t) - \big(R_{\varphi,m}f\big)(t) \big| \leq \big(E_1(m) + E_2(m)\big)\,\| f \|_{L^2(\R)}
\end{align*}
can be shown for the interval~\mbox{$[n,\,n+1]$} with arbitrary~\mbox{$n \in \Z$}. On the open interval~\mbox{$(n,\,n+1)$}, we decompose the approximation error as
\begin{align*}
f(t+n) - \big(R_{\varphi,m}f\big)(t)(t+n) = e_1(t+n) + e_{2,n}(t)\,, \quad t\in (0,\,1)\,,
\end{align*}
with
\begin{align*}
e_1(t+n) &= f(t+n) - \sum_{k\in \Z} f(k)\,\sinc\big(t- (k-n)\big)\,\varphi\big(t- (k-n)\big)\\
&= f(t+n) - \sum_{\ell\in \Z} f(\ell +n)\,\sinc(t-\ell)\,\varphi(t - \ell)\,,\\
e_{2,n}(t) &\coloneqq \sum_{\ell\in \Z \setminus \{1-m,\ldots, m\}}  f(\ell +n)\,\sinc(t-\ell)\,\varphi(t - \ell)\,.
\end{align*}
As shown in steps (ii) and (iii), we have
\begin{align*}
\|e_1(\cdot + n)\|_{C_0(\R)} &= \|e_1\|_{C_0(\R)}\,,\\[1ex]
|e_{2,n}(t)| &\leq E_2(m)\,\|f\|_{L^2(\R)}\,, \quad t \in (0,\,1)\,.
\end{align*}
Furthermore, by the interpolation property~\eqref{eq:interpolprop} of~\mbox{$R_{\varphi,m}f$}, we have~\mbox{$e_{2,n}(0) = e_{2,n}(1) = 0$} for each~\mbox{$n \in \Z$} and thus
\begin{align*}
\max_{t\in [n,n+1]} \big|e_{2,n}(t)\big| \leq E_2(m)\,\| f\|_{L^2(\R)}\,.
\end{align*}
Hence, it follows that
\begin{align*}
\max_{t \in [n,n+1]} \big| f(t) - \big(R_{\varphi,m}f\big)(t)\big| &\leq \|e_1\|_{C_0(\R)} + \max_{t\in [n,n+1]} \big|e_{2,n}(t)\big|\\
&\leq \big( E_1(m) + E_2(m)\big)\, \| f\|_{L^2(\R)}\,,
\end{align*}
which completes the proof. \hfill\qedsymbol

\section{Regularization with the Gaussian function \label{sec:Gauss}}

In this section we consider the Gaussian function~\eqref{eq:Gaussfunction} with variance~\mbox{$\sigma^2 > 0$}, analogous to~\cite[Theorem~4.1]{KiPoTa22}. In order to achieve fast convergence of the Gaussian regularized Shannon sampling formula, we \new{also study} the choice of this variance~$\sigma^2$.

\begin{Theorem}
\label{Theorem:approxerrorGauss}
Assume that~\mbox{$f\in L^2(\R) \cap C(\R)$} is bandlimited with bandwidth~\mbox{$\delta \in (0,\,\pi)$}. Further let~\mbox{$\varphi_{\mathrm{Gauss}}$} be the Gaussian function~\eqref{eq:Gaussfunction}
with variance~\mbox{$\sigma^2 = \frac{m}{\pi - \delta}$} and let~\mbox{$m \in {\mathbb N}\setminus \{1\}$} be given. \\
Then the Gaussian regularized Shannon sampling formula satisfies the error estimate
\begin{align}
\label{eq:err_Gauss}
\big\| f - R_{\mathrm{Gauss},m}f \big\|_{C_0(\R)} \le \frac{2\sqrt{2}}{\sqrt{\pi m\,(\pi -\delta)}}\, \e^{-m\,(\pi-\delta)/2}\,\| f\|_{L^2(\R)}\,.
\end{align}
\end{Theorem}

\emph{Proof}. (i) At first, we estimate the regularization error constant~\eqref{eq:E1} for the Gaussian function~\eqref{eq:Gaussfunction}.
Since the Fourier transform of~\mbox{$\varphi_{\mathrm{Gauss}}$} reads as
\begin{align*}
{\hat \varphi}_{\mathrm{Gauss}}(\omega) = \frac{1}{\sqrt{2\pi}}\, \int_{\R} \varphi_{\mathrm{Gauss}}(t)\, \e^{- \i\,t \omega}\,\mathrm{d}t = \sigma\,\e^{-\omega^2\sigma^2/2}\,, \quad \omega\in \R\,,
\end{align*}
cf.~\cite[Example~2.6]{PPST23}, we have
\begin{align*}
E_1(m) = \max_{\omega \in [-\delta,\delta]} \bigg| 1 - \frac{\sigma}{\sqrt{2\pi}}\,\int_{\omega - \pi}^{\omega + \pi} \e^{-\tau^2 \sigma^2/2}\,\mathrm{d}\tau\bigg|\,.
\end{align*}
Substituting~\mbox{$s = \tau \sigma/\sqrt 2$} and using the integral~\mbox{$\int_{\R} \e^{-s^2}\,\mathrm{d}s = \sqrt \pi$},
we obtain for~\mbox{$\omega \in [-\delta,\,\delta]$} with~\mbox{$\delta \in (0,\,\pi)$} that
\begin{align*}
\Delta_{\mathrm{Gauss}}(\omega) &\coloneqq 1 - \frac{1}{\sqrt \pi}\, \int_{(\omega - \pi)\sigma/\sqrt 2}^{(\omega + \pi)\sigma/\sqrt 2} \e^{-s^2}\,\mathrm{d}s\\
&= \frac{1}{\sqrt \pi}\,\bigg( \int_{\R} \e^{-s^2}\,\mathrm{d}s  -  \int_{(\omega - \pi)\sigma/\sqrt 2}^{(\omega + \pi)\sigma/\sqrt 2} \e^{-s^2}\,\mathrm{d}s\bigg)\\
&= \frac{1}{\sqrt \pi}\, \bigg( \int_{-\infty}^{(\omega - \pi)\sigma/\sqrt 2} \e^{-s^2}\,\mathrm{d}s + \int_{(\omega + \pi)\sigma/\sqrt 2}^{\infty} \e^{-s^2}\,\mathrm{d}s \bigg)\\
&= \frac{1}{\sqrt \pi}\, \bigg( \int_{(\pi - \omega)\sigma/\sqrt 2}^{\infty} \e^{-s^2}\,\mathrm{d}s + \int_{(\omega + \pi)\sigma/\sqrt 2}^{\infty} \e^{-s^2}\,\mathrm{d}s \bigg)\,.
\end{align*}
Since~\mbox{$\Delta_{\mathrm{Gauss}}$} is even, we consider only the case~\mbox{$\omega \in [0,\,\delta]$}. Applying the inequality
\begin{align*}
\int_a^{\infty}  \e^{-s^2}\,\mathrm{d}s = \int_0^{\infty} \e^{-(t+a)^2}\,\mathrm{d}t \leq \e^{-a^2}\,\int_0^{\infty} \e^{-2at}\,\mathrm{d}t = \frac{1}{2a}\,\e^{-a^2} \,,
\quad a > 0\,,
\end{align*}
we obtain
\begin{align*}
0 \leq \Delta_{\mathrm{Gauss}}(\omega) &\leq \frac{1}{\sqrt{2\pi}}\,\bigg(\frac{\e^{-(\pi-\omega)^2\sigma^2/2}}{(\pi-\omega)\,\sigma} + \frac{\e^{-(\pi+\omega)^2\sigma^2/2}}{(\pi+\omega)\,\sigma}\bigg) %\\ &
\leq \sqrt{\frac{2}{\pi}}\,\frac{\e^{-(\pi-\omega)^2\sigma^2/2}}{(\pi-\omega)\,\sigma} \,,
\quad \omega \in [0,\,\delta]\,.
\end{align*}
Consequently, we have for all~\mbox{$\omega \in [-\delta,\,\delta]$} that
\begin{align*}
0 \leq \Delta_{\mathrm{Gauss}}(\omega) \leq \sqrt{\frac{2}{\pi}}\,\frac{\e^{-(\pi-|\omega|)^2\sigma^2/2}}{(\pi-|\omega|)\,\sigma}
\end{align*}
and hence
\begin{align}
\label{eq:GaussE1}
E_1(m) \leq \sqrt{\frac{2}{\pi}}\,\frac{\e^{-(\pi-\delta)^2\sigma^2/2}}{(\pi-\delta)\,\sigma}\,.
\end{align}
(ii) Now we examine the truncation error constant~\eqref{eq:E2} for the Gaussian function~\eqref{eq:Gaussfunction}.
By~\mbox{$\varphi_{\mathrm{Gauss}}^2(m) = \e^{-m^2/\sigma^2}$} and the inequality
\begin{align*}
\int_m^{\infty} \varphi_{\mathrm{Gauss}}^2(t)\, \mathrm{d}t  = \sigma\, \int_{m/\sigma}^{\infty} \e^{-s^2}\,\mathrm{d}s \leq \frac{\sigma^2}{2m}\, \e^{-m^2/\sigma^2}
\end{align*}
we obtain
\begin{align}
\label{eq:GaussE2}
E_2(m) &\leq \frac{\sqrt 2}{\pi m}\,\sqrt{\e^{-m^2/\sigma^2} + \frac{\sigma^2}{2m}\, \e^{-m^2/\sigma^2}}
=
\frac{\sqrt 2}{\pi m}\,\sqrt{1 + \frac{\sigma^2}{2m}}\,{\e}^{-m^2/(2\sigma^2)} \,.
\end{align}
(iii) Finally, we \new{choose} the variance~\mbox{$\sigma^2$} of the Gaussian function~\eqref{eq:Gaussfunction} \new{such that}~\mbox{$E_1(m)$} and~\mbox{$E_2(m)$} possess the same exponential decay with respect to~$m$.
From~\eqref{eq:GaussE1} and~\eqref{eq:GaussE2} it follows that
\begin{align}
\label{eq:optvariance}
\sigma^2 \coloneqq \frac{m}{\pi - \delta}\,.
\end{align}
\new{This yields the estimates}
\begin{align*}
E_1(m) &\leq \sqrt{\frac{2}{\pi}}\, \frac{1}{\sqrt{m\,(\pi - \delta)}} \,\e^{-m\,(\pi - \delta)/2}\,,\\
E_2(m) &\leq \frac{\sqrt 2}{\pi m}\,\sqrt{1 + \frac{1}{2\,(\pi - \delta)}}\,\e^{-m\,(\pi - \delta)/2}\,.
\end{align*}
Note that since~\mbox{$m\in \N\setminus \{1\}$} and~\mbox{$ \delta \in (0,\,\pi)$}, we have
\begin{align*}
\left(\frac{\sqrt{2}}{\sqrt{\pi m\,(\pi - \delta)}}\right)^{-1} \cdot \frac{\sqrt 2}{\pi m}\,\sqrt{1 + \frac{1}{2\,(\pi - \delta)}}
&=
\sqrt{\frac{2(\pi-\delta)+1}{2\pi m}}
\leq
\sqrt{\frac{2\pi+1}{4\pi}}  < 1
\end{align*}
and therefore
\begin{align*}
E_2(m) \leq \frac{\sqrt{2}}{\sqrt{\pi m\,(\pi - \delta)}} \,\e^{-m\,(\pi - \delta)/2} \,.
\end{align*}
Thus, the Gaussian regularized Shannon sampling formula with the variance~\eqref{eq:optvariance} fulfills the error estimate~\eqref{eq:err_Gauss}.
This completes the proof. \hfill \qedsymbol
\medskip

Note that already in~\cite[Theorem~4.1]{KiPoTa22} bounds on the approximation error of the Shannon sampling formula~\eqref{eq:regShannonformula} were shown for the Gaussian function~\eqref{eq:Gaussfunction} with suitably chosen variance~\mbox{$\sigma^2$}, which is basically the same as the one in Theorem~\ref{Theorem:approxerrorGauss}, only looking slightly different due to the different setting considered in~\cite{KiPoTa22}.

\begin{Remark}
\label{Remark:opt_Gauss}	
\new{Inspired by~\cite{CZ19} one could define a weak form of optimality of the Gaussian regularized Shannon sampling formula by saying that the variance~\mbox{$\sigma^2$} of the Gaussian function~\eqref{eq:Gaussfunction} is \emph{optimal}, if~\mbox{$E_1(m)$} and~\mbox{$E_2(m)$} possess the same exponential decay with respect to~$m$.
Hence, Theorem~\ref{Theorem:approxerrorGauss} shows that the choice~\eqref{eq:optvariance} is optimal for the Shannon sampling formula~\eqref{eq:regShannonformula} with the Gaussian function~\eqref{eq:Gaussfunction} in this weak sense.}
We remark that in~\cite{CZ19} a \new{slightly} different optimal variance~\mbox{$\sigma^2 = \frac{m-1}{\pi - \delta}$} is presented for the Gaussian regularizer~\eqref{eq:Gaussfunction}, while \new{also considering} a slightly different truncation than in~\eqref{eq:regShannonformula}.
Nevertheless, both results, Theorem~\ref{Theorem:approxerrorGauss} and~\cite[Theorem~1.1]{CZ19}, possess the same asymptotic behavior.

Additionally, it should be noted that in~\cite{CZ19} the approximation error is estimated only up to an unknown constant, while our error estimate of the Gaussian regularized Shannon sampling formula contains relatively small explicit constants, which is more favorable for practical applications.
Moreover, we estimate the approximation error differently by splitting it into the regularization error~\eqref{eq:e1} and the truncation error~\eqref{eq:e20}, which seems more intuitive than the rather artificial analysis presented in~\cite[Theorem~1.1]{CZ19}.
\end{Remark}

\new{This definition of weak optimality has led to the following open question of optimality of the variance~\eqref{eq:optvariance}, which could so far only be observed numerically.
\begin{Conjecture}
\label{conj:opt_Gauss}
The parameter~\eqref{eq:optvariance} is the optimal variance for the Shannon sampling formula~\eqref{eq:regShannonformula} with the Gaussian function~\eqref{eq:Gaussfunction} not only in the weak sense of Remark~\ref{Remark:opt_Gauss}, but also guarantees the maximum decay rate of the uniform approximation error~\eqref{eq:err_approx}.
\end{Conjecture}}

\begin{Example}
\label{ex:opt_params_Gauss}
\new{In order to present numerical evidence for the optimality of the variance~\eqref{eq:optvariance} of the Gaussian regularized Shannon sampling formula stated in Conjecture~\ref{conj:opt_Gauss}} we consider the regularized Shannon sampling formula~\eqref{eq:regShannonformula} with the Gaussian function~\mbox{$\varphi_{\mathrm{Gauss}}$} in~\eqref{eq:Gaussfunction} \new{for a given bandlimited function~\mbox{$f \in L^2(\mathbb R) \cap C(\mathbb R)$} with bandwidth~\mbox{$\delta \in (0,\,\pi]$}} and \new{estimate} the corresponding approximation error
\begin{align}
\label{eq:maxerr}
\max_{t\in [-1,\, 1]} \big| f(t) - \big(R_{\varphi,m}f\big)(t)\big| \,,
\end{align}
cf.~\eqref{eq:err_approx}, \new{numerically}.
The error~\eqref{eq:maxerr} shall here be approximated by evaluating a given function~$f$ and its approximation~\mbox{$R_{\varphi,m}f$} at equidistant points~\mbox{$t_s\in[-1,\,1]$}, \mbox{$s=1,\dots,S$}, with~\mbox{$S=10^5$}.
Note that by the definition of the regularized Shannon sampling formula~\eqref{eq:regShannonformula} we have
\begin{align*}
\big(R_{\varphi,m}f\big)(t) = \sum_{k=-m-1}^{m+1} f(k)\, \sinc(t-k)\,\varphi(t-k) \,,
\quad t\in[-1,\,1] \,.
\end{align*}
Analogous to~\cite[Section IV, C]{Ob90} we study the bandlimited function
\begin{align}
\label{eq:testfunction}
f(t) = \frac{2\delta}{\sqrt{5\pi\delta+4\pi\sin\delta}} \left[\sinc\bigg(\frac{\delta t}{\pi}\bigg) + \frac 12 \,\sinc\bigg(\frac{\delta (t-1)}{\pi}\bigg)\right] \,,
\quad t\in\R \,,
\end{align}
with~\mbox{$\|f\|_{L^2(\R)}=1$}, for several bandwidth parameters~\mbox{$\delta \in \big\{\frac{\pi}{4},\, \frac{\pi}{2},\,\frac{3\pi}{4}\big\}$}, i.\,e., several oversampling rates~\mbox{$\frac{\pi}{\delta}>1$}.
To compare with the variance~\mbox{$\sigma^2$} \new{stated} in~\eqref{eq:optvariance}, we choose the parameter of the Gaussian function~\eqref{eq:Gaussfunction} as~\new{\mbox{$\sigma=\alpha \cdot \frac{m}{\pi - \delta}$}} with~\new{\mbox{$\alpha\in\big\{\frac 12,1,2\big\}$}}.

The corresponding results for different truncation parameters~\mbox{$m\in\{2, 3, \ldots, 10\}$} are displayed in Figure~\ref{fig:error_Gauss}.
It can clearly be seen that both, an increase and a decrease of the variance in~\eqref{eq:optvariance}, cause worsened error decay rates with respect to~$m$.
Thus, the numerical results \new{give reason to believe} that the variance~\eqref{eq:optvariance} of Theorem~\ref{Theorem:approxerrorGauss} is \new{indeed} optimal \new{in terms of the uniform approximation error~\eqref{eq:err_approx}}, already for very small truncation parameters~\mbox{$m\in\N\setminus\{1\}$}.
\begin{figure}[ht]
	\centering
	\captionsetup[subfigure]{justification=centering}
	\begin{subfigure}[t]{0.32\textwidth}
		\includegraphics[width=\textwidth]{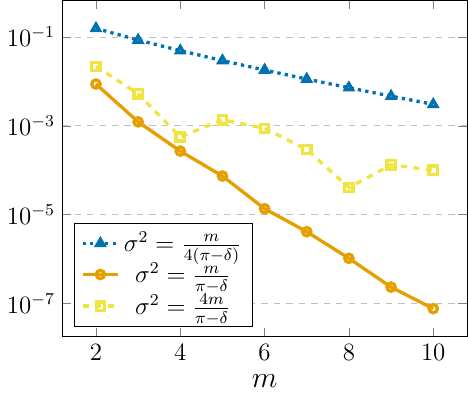}
		\caption{$\delta=\frac{\pi}{4}$}
	\end{subfigure}
	\begin{subfigure}[t]{0.32\textwidth}
		\includegraphics[width=\textwidth]{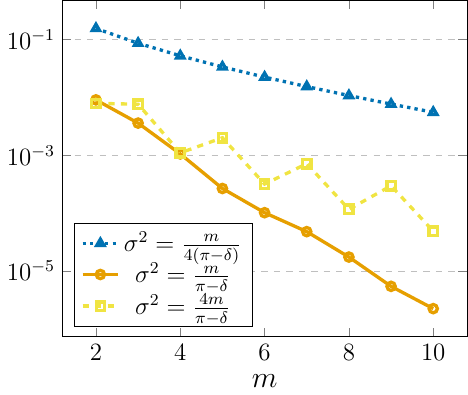}
		\caption{$\delta=\frac{\pi}{2}$}
	\end{subfigure}
	\begin{subfigure}[t]{0.32\textwidth}
		\includegraphics[width=\textwidth]{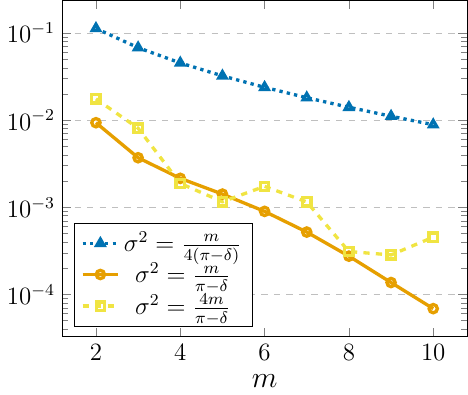}
		\caption{$\delta=\frac{3\pi}{4}$}
	\end{subfigure}
	\caption{Maximum approximation error~\eqref{eq:maxerr} using the Gaussian function~\mbox{$\varphi_{\mathrm{Gauss}}$} in~\eqref{eq:Gaussfunction} with different variances~\mbox{$\sigma^2\in\big\{\frac{m}{4(\pi-\delta)},\frac{m}{\pi-\delta},\frac{4m}{\pi-\delta}\big\}$}, for the bandlimited function~\eqref{eq:testfunction} with bandwidths~\mbox{$\delta \in \big\{\frac{\pi}{4},\frac{\pi}{2},\frac{3\pi}{4}\big\}$} and truncation parameters~\mbox{$m\in\{2, 3, \ldots, 10\}$}.
	\label{fig:error_Gauss}}
\end{figure}
\end{Example}

\begin{Remark}
\new{Note that~\cite{QO05, QC06} suggested the \emph{modified Gaussian function}
\begin{align}
\label{eq:modGaussfunction}
\varphi_{\mathrm{modGauss}}(t) \coloneqq \e^{-t^2/(2 \sigma^2)}\,\cos(\lambda t)\,, \quad t \in \mathbb R\,,
\end{align}
with the parameters~\mbox{$\sigma^2 > 0$} and~\mbox{$\lambda \geq 0$} as an improvement to the Gaussian function~\eqref{eq:Gaussfunction}.
By the same techniques as in Theorem~\ref{Theorem:approxerrorGauss}, however, one can determine that the optimal variance of~\eqref{eq:modGaussfunction} in the weak sense of Remark~\ref{Remark:opt_Gauss} is given by~\mbox{$\sigma^2 = \frac{m}{\pi - \lambda - \delta}$}, \mbox{$0 \leq \lambda < \pi - \delta$}, with the corresponding error estimate
\begin{align*}
\big\| f - R_{\mathrm{modGauss},m}f \big\|_{C_0(\R)} \le \frac{2\sqrt{2}}{\sqrt{\pi m\,(\pi - \lambda - \delta)}}\, \e^{-m\,(\pi-\lambda-\delta)/2}\,\| f\|_{L^2(\R)}\,.
\end{align*}
This shows that the approximation error of the regularized Shannon sampling formula with the modified Gaussian function~\eqref{eq:modGaussfunction} has the best exponential decay in the case~\mbox{$\lambda = 0$}, therefore proving that the Gaussian function~\mbox{$\varphi_{\mathrm{Gauss}}$} in~\eqref{eq:Gaussfunction} is much more favorable than the modified Gaussian function~\mbox{$\varphi_{\mathrm{modGauss}}$} in~\eqref{eq:modGaussfunction}.}
\end{Remark}

\section{Regularization with the sinh-type window function \label{sec:sinh}}

In this section, we consider the $\sinh$-type window function~\eqref{eq:varphisinh} with shape parameter~\mbox{$\beta > 0$}, analogous to~\cite[Theorem~6.1]{KiPoTa22} and~\cite[Theorem~4.2]{KiPoTa24}.
\new{We especially focus on addressing minor gaps in~\cite{KiPoTa22,KiPoTa24} by rigorously proving assumptions up to now based solely on numerical tests.}
Moreover, we demonstrate that the exponential decay with respect to the truncation parameter~\mbox{$m \in \N \setminus \{1\}$} is \new{twice as fast} for the uniform approximation error~\mbox{$\big\| f - R_{\sinh,m}f \big\|_{C_0(\R)}$} \new{as} for the approximation error~\mbox{$\big\| f - R_{\mathrm{Gauss},m}f \big\|_{C_0(\R)}$} in Theorem~\ref{Theorem:approxerrorGauss}.
\new{To this end, we firstly formulate the following lemma.
\begin{Lemma}
\label{Lemma:estimate_integralJ1}
For all~\mbox{$W > 1$} and~\mbox{$\beta>0$} we have
\begin{align}
	\label{eq:integralJ1_est}
	\bigg| \int_1^W \frac{J_1(\beta\,\sqrt{\nu^2-1})}{\sqrt{\nu^2 - 1}}\,{\mathrm d}\nu \,\bigg|
	\leq
	\frac{1 - {\mathrm e}^{-\beta}}{\beta} + \frac{{\sqrt 2}\,\pi}{\sqrt \beta}\,,
\end{align}
where~$J_1$ denotes the Bessel function of first order.
\end{Lemma}}

\new{\emph{Proof.} Substituting~\mbox{$\nu = \cosh t$} in~\eqref{eq:integralJ1_est}, we obtain
\begin{align*}
	\int_1^W \frac{J_1(\beta\,\sqrt{\nu^2-1})}{\sqrt{\nu^2 - 1}}\,{\mathrm d}\nu
	= \int_0^w J_1(\beta\,\sinh t)\,{\mathrm d}t
\end{align*}
with~\mbox{$w = \mathrm{arcosh} W > 0$}.
Note that it is known by~\cite[6.645--1]{GR80} that
\begin{align*}
	\int_0^{\infty} J_1(\beta\,\sinh t)\,{\mathrm d}t
	= I_{1/2}\Big(\frac{\beta}{2}\Big)\, K_{1/2}\Big(\frac{\beta}{2}\Big)
	= \frac{1 - {\mathrm e}^{-\beta}}{\beta} > 0\,,
\end{align*}
where~\mbox{$I_{1/2}$} and~\mbox{$K_{1/2}$} denote the modified Bessel functions of half order (see~\cite[10.2.13, 10.2.14, and~10.2.17]{abst}.
The additional substitution~\mbox{$s = \beta\,\sinh t$} yields
\begin{align}
	\label{eq:integralJ1_infty}
	\int_0^{\infty} J_1(\beta\,\sinh t)\,{\mathrm d}t
	&= \int_0^{\infty} \frac{J_1(s)}{\sqrt{s^2 + \beta^2}}\,{\mathrm d}s
	= \frac{1 - {\mathrm e}^{-\beta}}{\beta}\,, \\ \notag
	\int_0^w J_1(\beta\,\sinh t)\,{\mathrm d}t
	&= \int_0^u \frac{J_1(s)}{\sqrt{s^2 + \beta^2}}\,{\mathrm d}s\,,
\end{align}
where~\mbox{$u = \beta\,\sinh w > 0$}.}

\new{Let~\mbox{$j_k$}, \mbox{$k \in \mathbb N$}, denote the positive zeros of~\mbox{$J_1$}. 
Note that~\mbox{$j_k$}, \mbox{$k=1,\dots,40$}, are tabulated in~\cite[p.~748]{Wa44} and
that by~\mbox{$(-1)^k\,J_1^{\prime}(j_k) > 0$}, see~\cite[9.1.27 and~9.5.2]{abst}, these zeros are simple.
Using these zeros the integral~\eqref{eq:integralJ1_infty} can be represented as
\begin{align}
\label{eq:integralJ1_series_representation}
\frac{1 - {\mathrm e}^{-\beta}}{\beta}
=
\int_0^{\infty} \frac{J_1(s)}{\sqrt{s^2 + \beta^2}}\,{\mathrm d}s
= \bigg(\int_0^{j_1} + \int_{j_1}^{j_2} + \int_{j_2}^{j_3} +\,\ldots \bigg)\frac{J_1(s)}{\sqrt{s^2 + \beta^2}}\,{\mathrm d}s\,.
\end{align}
The integrand~\mbox{$J_1(s)\,(s^2 + \beta^2)^{-1/2}$}, \mbox{$s \in [0,\infty)$}, is an oscillating function which tends to zero for~\mbox{$s \to \infty$}, since we have
\begin{align}
\label{eq:est_J1}
\frac{|J_1(s)|}{\sqrt{s^2 + \beta^2}} \leq \frac{1}{{\sqrt s}\,\sqrt{s^2 + \beta^2}}\,,
\quad s \in [j_1,\, \infty)\,,
\end{align}
which can be shown by the same technique as in~\cite[Lemma~6]{PT21b}.
Thus, the right-hand side of~\eqref{eq:integralJ1_series_representation} is a convergent alternating series.}

\new{Now we consider the function
\begin{align}
	\label{eq:integral_J1}
	\mathcal I(u) \coloneqq \int_0^u \frac{J_1(s)}{\sqrt{s^2 + \beta^2}}\,{\mathrm d}s\,,
	\quad u \in [0,\,\infty)\,.
\end{align}
We immediately recognize that~\mbox{$\mathcal I(u)$} has relative maxima at~\mbox{$j_{2n+1}$}, \mbox{$n\in\N$}, and relative minima at~\mbox{$j_{2n}$}, \mbox{$n\in\N$}, since we have
\begin{align*}
	\mathcal I^{\prime}(j_k) = \frac{J_1(j_k)}{\sqrt{j_k^2 + \beta^2}}
	= 0\,, \quad (-1)^k\,\mathcal I^{\prime \prime}(j_k)
	= \frac{(-1)^k\,J_1^{\prime}(j_k)}{\sqrt{j_k^2 + \beta^2}} > 0 \,,
	\quad k\in\N\,.
\end{align*}
In other words, by the oscillatory behavior of the Bessel function~\mbox{$J_1$}, the function~\mbox{$\mathcal I(u)$} increases from~\mbox{$\mathcal I(0) = 0$} to~\mbox{$\mathcal I(j_1)$}, then decreases from~\mbox{$\mathcal I(j_1)$} to~\mbox{$\mathcal I(j_2)$}, increases again from~\mbox{$\mathcal I(j_2)$} to~\mbox{$\mathcal I(j_3)$}, and so on.
By~\eqref{eq:integralJ1_infty} and~\eqref{eq:integral_J1} it is also easy to see that
\begin{align*}
	\frac{1 - {\mathrm e}^{-\beta}}{\beta} - \mathcal I(j_k)
	=
	\bigg(\int_0^{\infty}-\int_0^{j_k}\bigg) \frac{J_1(s)}{\sqrt{s^2 + \beta^2}}\,{\mathrm d}s
	=
	\int_{j_k}^{\infty} \frac{J_1(s)}{\sqrt{s^2 + \beta^2}} \,{\mathrm d}s \,,
	\quad k\in\N\,.
\end{align*}
Thus, combined with the estimate~\eqref{eq:est_J1} this yields 
\begin{align*}
	\bigg| \frac{1 - {\mathrm e}^{-\beta}}{\beta} - \mathcal I(j_k)\bigg|
	\leq \int_{j_k}^{\infty} \frac{|J_1(s)|}{\sqrt{s^2 + \beta^2}} \,{\mathrm d}s
	\leq \int_{j_k}^{\infty} \frac{1}{{\sqrt s}\,\sqrt{s^2 + \beta^2}} \,{\mathrm d}s \,,
	\quad k\in\N\,.
\end{align*}
In addition, from the equivalence relation~\mbox{$\|\boldsymbol v\|_1 \leq \sqrt{n} \,\|\boldsymbol v\|_2$}, \mbox{$\boldsymbol v \in \R^n$}, of the vector norms it follows that~\mbox{$s + \beta \leq \sqrt 2 \,\sqrt{s^2 + \beta^2}$}, \mbox{$s, \beta>0$}, and therefore
\begin{align*}
	\frac{1}{\sqrt{s^2 + \beta^2}} \leq \frac{\sqrt 2}{s + \beta}\,,
	\quad s, \beta>0\,.
\end{align*}
Hence, we obtain
\begin{align*}
	\bigg| \frac{1 - {\mathrm e}^{-\beta}}{\beta} - \mathcal I(j_k)\bigg| \leq \int_{j_k}^{\infty}\frac{\sqrt{2}}{{\sqrt s}\,(s + \beta)} \,{\mathrm d}s \,,
	\quad k\in\N\,.
\end{align*}
Since the antiderivative of~\mbox{$s^{-1/2}\,(s + \beta)^{-1}$} reads as~\mbox{$\frac{2}{\sqrt \beta}\,\arctan \sqrt{\frac{s}{\beta}}$} and~\mbox{$\arctan y \leq \frac{\pi}{2}$}, \mbox{$y\in\R$}, we obtain the estimate
\begin{align*}
	\bigg| \frac{1 - {\mathrm e}^{-\beta}}{\beta} - \mathcal I(j_k)\bigg|
	&\leq \frac{2\,\sqrt 2}{\sqrt \beta}\,\Bigg(\frac{\pi}{2} - \arctan \sqrt{\frac{j_k}{\beta}}\,\Bigg) %\notag \\
	= \frac{2\,\sqrt 2}{\sqrt \beta}\,\arctan \sqrt{\frac{\beta}{j_k}}
	\leq \frac{\sqrt 2 \,\pi}{\sqrt \beta}\,.
\end{align*}
Consequently, as this estimate is valid for all relative extreme values~\mbox{$\mathcal I(j_k)$}, \mbox{$k \in \mathbb N$}, we also obtain
\begin{align*}
%\label{eq:1}
	\bigg| \frac{1 - {\mathrm e}^{-\beta}}{\beta} - \mathcal I(u) \bigg| \leq \frac{{\sqrt 2}\,\pi}{\sqrt \beta}
\end{align*}
for all~\mbox{$u>0$}, which immediately implies
\begin{align*}
	\frac{1 - {\mathrm e}^{-\beta}}{\beta} - \frac{{\sqrt 2}\,\pi}{\sqrt \beta}
	\leq \mathcal I(u)
	\leq \frac{1 - {\mathrm e}^{-\beta}}{\beta} + \frac{{\sqrt 2}\,\pi}{\sqrt \beta}\,.
\end{align*}
Since~\mbox{$\frac{1 - {\mathrm e}^{-\beta}}{\beta} > 0$} and~\mbox{$\frac{{\sqrt 2}\,\pi}{\sqrt \beta} > 0$} for~\mbox{$\beta > 0$} this yields
\begin{align*}
	- \frac{1 - {\mathrm e}^{-\beta}}{\beta} - \frac{{\sqrt 2}\,\pi}{\sqrt \beta}
	\leq \frac{1 - {\mathrm e}^{-\beta}}{\beta} - \frac{{\sqrt 2}\,\pi}{\sqrt \beta}
	\leq \mathcal I(u)
	\leq \frac{1 - {\mathrm e}^{-\beta}}{\beta} + \frac{{\sqrt 2}\,\pi}{\sqrt \beta}\,,
\end{align*}
and thereby the assertion~\eqref{eq:integralJ1_est}. 
This completes the proof. \hfill \qedsymbol
\medskip}

\begin{Theorem}
\label{Theorem:approxerrorsinh}
Let~\mbox{$m \in {\mathbb N}\setminus \{1\}$} be given.
Assume that~\mbox{$f\in L^2(\R) \cap C(\R)$} is bandlimited with bandwidth~\mbox{$\delta \in \big( 0,\,\new{\frac{m-1}{m}\,\pi} \big]$}. Further let~\mbox{$\varphi_{\sinh}$} be the $\sinh$-type window function~\eqref{eq:varphisinh}
with shape parameter~\mbox{$\beta = m\,(\pi - \delta)$}. \\
Then the $\sinh$-type regularized Shannon sampling formula satisfies the error estimate
\begin{align}
\label{eq:err_sinh}
\big\| f - R_{\sinh,m}f \big\|_{C_0(\R)} \le \new{\big(4+9\,\sqrt{m\,(\pi-\delta)}\,\big)} \,\e^{-m\,(\pi-\delta)}\,\| f\|_{L^2(\R)}\,.
\end{align}
\end{Theorem}

\emph{Proof}. (i) Since~\mbox{$\varphi_{\sinh}$} in~\eqref{eq:varphisinh} is compactly supported on~\mbox{$[-m,\,m]$} and~\mbox{$\varphi_{\sinh}(m) = 0$}, we have~\mbox{$E_2(m) = 0$}. Thus, according to Theorem~\ref{Theorem:approxerror}, the approximation
error can be estimated by
\begin{align*}
\big\| f - R_{\sinh,m}f \big\|_{C_0(\R)} \le \| f\|_{L^2(\R)}\, \max_{\omega \in [-\delta,\delta]} \big| \Delta_{\sinh}(\omega)\big|\,,
\end{align*}
where
\begin{align}
\label{eq:Deltasinh}
\Delta_{\sinh}(\omega) \coloneqq 1 - \frac{1}{\sqrt{2\pi}}\,\int_{\omega - \pi}^{\omega+\pi} {\hat \varphi}_{\sinh}(\tau)\,\mathrm{d}\tau\,, \quad \omega \in [-\delta,\, \delta]\,.
\end{align}
Following~\cite[p.~38, 7.58]{Ob90}, the Fourier transform of~\eqref{eq:varphisinh} has the form
\begin{align}
\label{eq:hatvarphisinh}
{\hat \varphi}_{\sinh}(\tau) = \frac{m\,\sqrt \pi}{\sqrt{2}\,\sinh \beta}\cdot \begin{cases} (1-\nu^2)^{-1/2}\,I_1\big(\beta\,\sqrt{1-\nu^2}\big) &\colon |\nu| < 1\,, \\[1ex]
 (\nu^2-1)^{-1/2}\,J_1\big(\beta\,\sqrt{\nu^2-1}\big) &\colon |\nu| > 1\,,
\end{cases}
\end{align}
with the scaled frequency~\mbox{$\nu =\frac{m}{\beta}\,\tau$}, where~$J_1$ denotes the Bessel function and~$I_1$ the modified Bessel function of first order.
Substituting~\mbox{$\tau = \frac{\beta}{m}\,\nu$} in the integral in~\eqref{eq:Deltasinh}, the function~\mbox{$\Delta_{\sinh}$} reads as
\begin{align}
\label{eq:Deltasinh_subst}
\Delta_{\sinh}(\omega) \coloneqq 1 - \frac{\beta}{\sqrt{2\pi}\,m}\,\int_{-\nu_1(-\omega)}^{\nu_1(\omega)} {\hat \varphi}_{\sinh}\big(\tfrac{\beta}{m}\,\nu\big)\,\mathrm{d}\nu\,, \quad \omega \in [-\delta,\, \delta]\,,
\end{align}
with the increasing linear function
\begin{align}
\label{eq:nu1}
	\nu_1(\omega) \coloneqq \frac{m}{\beta}\, (\omega + \pi)\,, \quad \omega \in [-\delta,\, \delta]\,.
\end{align}

(ii) Now we choose the shape parameter of~\eqref{eq:varphisinh} in the special form~\mbox{$\beta = m\,(\pi - \delta)$}. Thus, we have
\begin{align*}
1 = \nu_1(-\delta) \leq \nu_1(\omega) = \frac{\omega + \pi}{\pi - \delta} \leq \nu_1(\delta) = \frac{\pi + \delta}{\pi - \delta}\,, \quad \omega \in [- \delta,\, \delta]\,.
\end{align*}
In view of~\eqref{eq:hatvarphisinh} we split~\eqref{eq:Deltasinh_subst} in the form~\mbox{$\Delta_{\sinh}(\omega) = \Delta_{\sinh,1} - \Delta_{\sinh,2}(\omega)$} with
\begin{align*}
\Delta_{\sinh,1} &\coloneqq 1 - \frac{\beta}{\sinh \beta}\,\int_0^{1} \frac{I_1\big(\beta\,\sqrt{1 - \nu^2}\big)}{\sqrt{1 - \nu^2}}\,\mathrm{d}\nu\,, \\[1ex]
\Delta_{\sinh,2}(\omega) &\coloneqq \frac{\beta}{2\,\sinh \beta}\,\bigg(\int_1^{\nu_1(-\omega)} + \int_1^{\nu_1(\omega)}\bigg)\,\frac{J_1\big(\beta\,\sqrt{\nu^2-1}\big)}{\sqrt{\nu^2-1}}\,\mathrm{d}\nu\,.
\end{align*}
Using~\cite[6.681--11]{GR80} and~\cite[10.2.13]{abst}, we get
\begin{align*}
\int_0^{1} \frac{I_1\big(\beta\,\sqrt{1 - \nu^2}\big)}{\sqrt{1 - \nu^2}}\,\mathrm{d}\nu = \int_0^{\pi/2} I_1(\beta\,\cos \sigma)\, \mathrm{d}\sigma = \frac{\pi}{2}\,\bigg(I_{1/2}\bigg(\frac{\beta}{2}\bigg)\bigg)^2 =
\frac{2}{\beta} \left(\sinh \frac{\beta}{2}\right)^2
\end{align*}
and hence
\begin{align}
\label{eq:shapecond1}
\new{0 < \;} \Delta_{\sinh,1} = 1 - \frac{2\,\big(\sinh \frac{\beta}{2}\big)^2}{\sinh \beta} = \frac{2\, \e^{-\beta}}{1 + \e^{-\beta}} \new{\;< 2\e^{-\beta}}\,.
\end{align}
\new{By Lemma~\ref{Lemma:estimate_integralJ1} we have}
\begin{align}
\label{eq:shapecond2}
\new{\big|\Delta_{\sinh,2}(\omega)\big|
\leq \frac{\beta}{\, \sinh \beta} \left( \frac{1 - {\mathrm e}^{-\beta}}{\beta} + \frac{{\sqrt 2}\,\pi}{\sqrt \beta} \right)
< \left( 2 + \sqrt{\beta} \,\frac{2\,\sqrt{2}\,\pi}{1-\e^{-2\pi}} \right) \!\e^{-\beta}}\,,
\quad \omega \in [- \delta,\, \delta]\,,
\end{align}
\new{since by assumption the parameters~\mbox{$m \in {\mathbb N}\setminus \{1\}$} and~\mbox{$\delta \in \big( 0,\,\frac{m-1}{m}\,\pi \big]$} are chosen such that \mbox{$\beta = m\,(\pi -\delta) \geq \pi$}.
Additionally, note that
\begin{align*}
\frac{2\,\sqrt{2}\,\pi}{1-\e^{-2\pi}} = 8.902390\ldots < 9
\end{align*}
holds.}
Thereby, the terms~\eqref{eq:shapecond1} and~\eqref{eq:shapecond2} have the same exponential decay~\mbox{$m\,(\pi - \delta)$} and~\eqref{eq:Deltasinh_subst} can be estimated by
\begin{align*}
\big| \Delta_{\sinh}(\omega)\big| = \new{\Delta_{\sinh,1} + \big|\Delta_{\sinh,2}(\omega)\big| \leq \big( 4+9\,\sqrt{\beta} \,\big) \,\e^{-\beta}}\,, \quad \omega \in [-\delta,\,\delta]\,.
\end{align*}
Thus, the $\sinh$-type regularized Shannon sampling formula with the chosen shape parameter~\mbox{$\beta = m\,(\pi -\delta)$} fulfills the error estimate~\eqref{eq:err_sinh}.
This completes the proof. \hfill \qedsymbol
\medskip

\new{Note} that already in~\cite[Theorem~6.1]{KiPoTa22} and~\cite[Theorem~4.2]{KiPoTa24} bounds on the approximation error of the Shannon sampling formula~\eqref{eq:regShannonformula} were shown for the \mbox{$\sinh$-type} window function~\eqref{eq:varphisinh} with suitably chosen shape parameter~$\beta$.
Although the respective parameters~$\beta$ look different than the one in Theorem~\ref{Theorem:approxerrorsinh}, they are basically the same, only adapted to the slightly different settings considered in~\cite{KiPoTa22, KiPoTa24}.
\new{However, it has to be pointed out that the error constant in Theorem~\ref{Theorem:approxerrorsinh} is somewhat worsened in comparison to our previous findings in~\cite{KiPoTa22, KiPoTa24} due to the fact that Lemma~\ref{Lemma:estimate_integralJ1} comprises a weaker version of the numerical assumption in~\cite[p.~25]{KiPoTa22}, but nevertheless closes the gap in this previous proof.}

\new{In addition, similar to Section~\ref{sec:Gauss}, the optimality of the shape parameter~\mbox{$\beta = m\,(\pi - \delta)$} for the $\sinh$-type window function~\eqref{eq:varphisinh} is still an open problem, which could so far only be observed numerically.
\begin{Conjecture}
\label{conj:opt_sinh}
The parameter~\mbox{$\beta = m\,(\pi - \delta)$} is optimal for the Shannon sampling formula~\eqref{eq:regShannonformula} with the $\sinh$-type window function~\eqref{eq:varphisinh}, as it guarantees the maximum decay rate of the uniform approximation error~\eqref{eq:err_approx}.
\end{Conjecture}}

\begin{Example}
\label{ex:opt_params_sinh}
Analogously as in Example~\ref{ex:opt_params_Gauss} we now \new{show numerical evidence for} the optimality of the shape parameter~\mbox{$\beta = m\,(\pi - \delta)$} of the $\sinh$-type regularized Shannon sampling formula \new{stated} in \new{Conjecture~\ref{conj:opt_sinh}}.
More precisely, for the bandlimited function~\eqref{eq:testfunction} with several bandwidth parameters~\mbox{$\delta \in \big\{\frac{\pi}{4},\, \frac{\pi}{2},\,\frac{3\pi}{4}\big\}$}, i.\,e., several oversampling rates~\mbox{$\frac{\pi}{\delta}>1$}, we consider the regularized Shannon sampling formula~\eqref{eq:regShannonformula} with the \mbox{$\sinh$-type} window function~\mbox{$\varphi_{\sinh}$} in~\eqref{eq:varphisinh}.
The corresponding approximation error~\eqref{eq:maxerr} shall again be approximated by evaluating the given function~$f$ and its approximation~\mbox{$R_{\varphi,m}f$} at equidistant points~\mbox{$t_s\in[-1,\,1]$}, \mbox{$s=1,\dots,S$}, with~\mbox{$S=10^5$}.
To compare with the parameter \new{in Theorem~\ref{Theorem:approxerrorsinh}}, we choose the shape parameter of the \mbox{$\sinh$-type} window function~\eqref{eq:varphisinh} as~\mbox{$\beta = \alpha\, m\,(\pi - \delta)$} with~\mbox{$\alpha\in\big\{\frac 12, 1, 2\big\}$}.

The outcomes for different truncation parameters~\mbox{$m\in\{2, 3, \ldots, 10\}$} are depicted in Figure~\ref{fig:error_sinh}.
\new{As supposed} it can clearly be seen that the choice of~\mbox{$\alpha \neq 1$} causes worsened error decay rates with respect to~$m$.
Thus, these numerical results \new{bolster the assertion} that the shape parameter~\mbox{$\beta = m\,(\pi - \delta)$} of Theorem~\ref{Theorem:approxerrorsinh} is \new{indeed} optimal \new{in terms of the uniform approximation error~\eqref{eq:err_approx}}, already for very small truncation parameters~\mbox{$m\in\N\setminus\{1\}$}.
\begin{figure}[ht]
	\centering
	\captionsetup[subfigure]{justification=centering}
	\begin{subfigure}[t]{0.32\textwidth}
		\includegraphics[width=\textwidth]{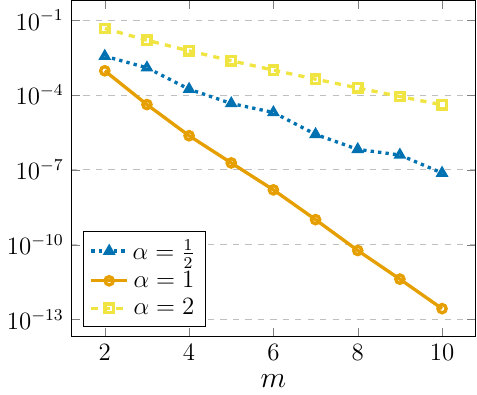}
		\caption{$\delta=\frac{\pi}{4}$}
	\end{subfigure}
	\begin{subfigure}[t]{0.32\textwidth}
		\includegraphics[width=\textwidth]{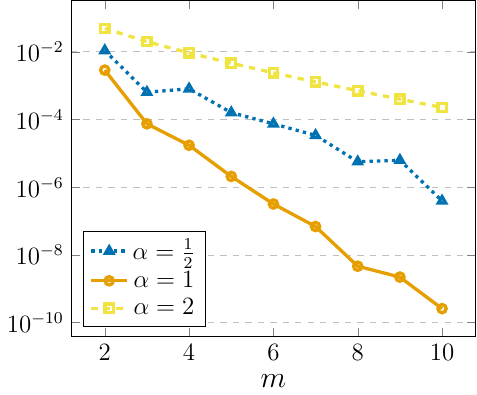}
		\caption{$\delta=\frac{\pi}{2}$}
	\end{subfigure}
	\begin{subfigure}[t]{0.32\textwidth}
		\includegraphics[width=\textwidth]{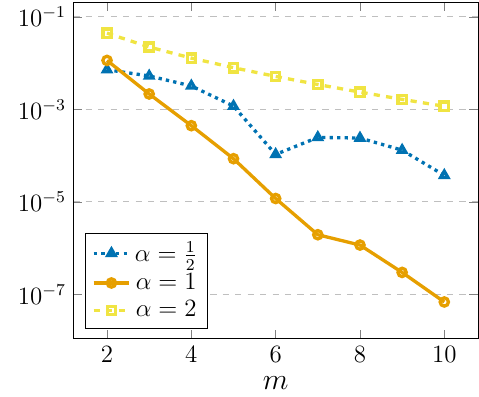}
		\caption{$\delta=\frac{3\pi}{4}$}
	\end{subfigure}
	\caption{Maximum approximation error~\eqref{eq:maxerr} using the \mbox{$\sinh$-type} window function~\mbox{$\varphi_{\mathrm{sinh}}$} in~\eqref{eq:varphisinh} with different shape parameters~\mbox{$\beta=\alpha\, m (\pi-\delta)$}, \mbox{$\alpha\in\big\{\frac 12,1,2\big\}$}, for the bandlimited function~\eqref{eq:testfunction} with bandwidths~\mbox{$\delta \in \big\{\frac{\pi}{4},\frac{\pi}{2},\frac{3\pi}{4}\big\}$} and truncation parameters~\mbox{$m\in\{2, 3, \ldots, 10\}$}.
	\label{fig:error_sinh}}
\end{figure}
\end{Example}

\section{Regularization with the continuous Kaiser--Bessel window function \label{sec:cKB}}

In this section, we consider the continuous Kaiser--Bessel window function~\eqref{eq:varphicKB} with shape parameter~\mbox{$\beta > 0$}, analogous to~\cite[Theorem~4.3]{KiPoTa24}.
\new{Once more, we put special emphasis on addressing minor gaps in~\cite{KiPoTa24} by rigorously proving assumptions that were previously based solely on numerical tests.}
Furthermore, we show that the exponential decay with respect to the truncation parameter~\mbox{$m \in \N \setminus \{1\}$} for the uniform approximation error~\mbox{$\big\| f - R_{\mathrm{cKB},m}f \big\|_{C_0(\R)}$} is \new{the same as for} the approximation error~\mbox{$\big\| f - R_{\sinh,m}f \big\|_{C_0(\R)}$} in Theorem~\ref{Theorem:approxerrorsinh}\new{, only with a slightly worse error constant.
To do so, we firstly establish the following lemmas.
\begin{Lemma}
\label{Lemma:abs_I0(beta)-L0(beta)}
For all~\mbox{$\beta\geq\pi$} we have
\begin{align}
\label{eq:abs_I0(beta)-L0(beta)}
\left|I_0(\beta) - {\mathbf L}_0(\beta) -1 + \frac{2}{\pi}\, \mathrm{Si}(\beta)\right| \leq \frac 12 \,,
\end{align}
where~\mbox{${\textbf L}_0$} denotes the \emph{modified Struve function}
\begin{align}
\label{eq:L0(x)}
{\textbf L}_0(x) \coloneqq \sum_{k=0}^{\infty} \frac{(x/2)^{2k+1}}{\big(\Gamma\big(k + \frac{3}{2}\big)\big)^2} = \frac{2x}{\pi}\, \sum_{k=0}^{\infty} \frac{x^{2k}}{\big((2k+1)!!\big)^2}\,, \quad x \in \mathbb R\,,
\end{align}
see~\cite[12.2.1]{abst}, and~\mbox{$\mathrm{Si}$} is the \emph{sine integral function}
\begin{align}
\label{eq:Si(x)}
\mathrm{Si}(x) \coloneqq \int_0^{x} \frac{\sin v}{v}\,{\mathrm d}v\,, \quad x \in \mathbb R\,.
\end{align}
\end{Lemma}
\emph{Proof}. 
By~\cite[Theorem~1]{BP14} the function~\mbox{$h(x) \coloneqq I_0(x) - {\mathbf L}_0(x)$} is completely monotonic on~\mbox{$[0,\,\infty)$}, i.\,e., it satisfies~\mbox{$h(x) \geq 0$} and~\mbox{$h^{\prime}(x) \leq 0$} for all~\mbox{$x\in[0,\,\infty)$}.
Thereby, we have
\begin{align}
\label{eq:ineq_hx}
0 \leq h(\beta) \leq h(\pi) \,, \quad \beta \geq \pi\,.
\end{align}
Due to the fact that
\begin{align*}
\mathrm{Si}^\prime(x) = \sinc\big(\tfrac{x}{\pi}\big)
\quad\text{and}\quad
\mathrm{Si}^{\prime\prime}(x) = \begin{cases}
\frac{x\cos x - \sin x}{x^2} & \colon x \neq 0\,, \\
0 & \colon x=0\,,
\end{cases}
\end{align*}
the sine integral function~\eqref{eq:Si(x)} has its local maxima at~\mbox{$(2k+1)\pi$}, \mbox{$k\in\Z\setminus\{0\}$}, and its local minima at~\mbox{$2k\pi$}, \mbox{$k\in\Z\setminus\{0\}$}.
Moreover, by the definition~\eqref{eq:Si(x)} the extremal points of the sine integral function become smaller in magnitude for~\mbox{$x\geq\pi$} when~\mbox{$x\to\infty$}, and~\mbox{$\lim_{x\to\infty}\mathrm{Si}(x) = \frac{\pi}{2}$}.
Thus, for all~\mbox{$\beta \geq \pi$} we have
\begin{align*}
\mathrm{Si}(2 \pi) \leq \mathrm{Si}(\beta) \leq \mathrm{Si}(\pi)\,,
\end{align*}
and consequently
\begin{align}
\label{eq:ineq_Si}
-1 + \frac{2}{\pi}\,\mathrm{Si}(2 \pi) \leq -1 + \frac{2}{\pi}\, \mathrm{Si}(\beta) \leq -1 + \frac{2}{\pi}\,\mathrm{Si}(\pi) \,, \quad \beta \geq \pi\,.
\end{align}
Combining~\eqref{eq:ineq_hx} and~\eqref{eq:ineq_Si}, we obtain the inequality
\begin{align*}
-0.097176... = -1 + \frac{2}{\pi}\,\mathrm{Si}(2 \pi) &\leq h(\beta) -1 + \frac{2}{\pi}\, \mathrm{Si}(\beta) \\ &\leq h(\pi) -1 + \frac{2}{\pi}\,\mathrm{Si}(\pi) = 0.400229\ldots
\end{align*}
for all~\mbox{$\beta \geq \pi$}, which finally gives~\eqref{eq:abs_I0(beta)-L0(beta)}.
This completes the proof.
\hfill\qedsymbol
\medskip
\begin{Lemma}
\label{Lemma:monotonicity}
The function
\begin{align}
\label{eq:func_g}
g(x) \coloneqq \frac{\e^x}{x\,(I_0(x)-1)}
\end{align}
is monotonously decreasing for~\mbox{$x\geq 1$}.
\end{Lemma}
\emph{Proof}.
Using~\mbox{$I_0^{\prime}(x) = I_1(x)$}, see~\cite[9.6.27]{abst}, the derivative of the function~\eqref{eq:func_g} is given by
\begin{align*}
g^\prime(x) = \frac{\e^x\,(-x-x\,I_1(x)+(x-1)\,I_0(x)+1)}{x^2\,(I_0(x)-1)^2} \,.
\end{align*}
To prove that~\eqref{eq:func_g} is monotonously decreasing for~\mbox{$x\geq 1$}, we need to show that~\mbox{$g^\prime(x)<0$}, \mbox{$x\geq 1$}, or rather
\begin{align*}
\tilde g(x) \coloneqq -x-x\,I_1(x)+(x-1)\,I_0(x)+1 < 0, \quad x \geq 1 \,.
\end{align*}
Note that
\begin{align*}
\tilde g(1) = -1-I_1(1)+0+1 =-I_1(1) = -0.565159\ldots < 0 \,.
\end{align*}
Thus, by showing~\mbox{$\tilde g^\prime(x) \leq 0$}, \mbox{$x\geq 1$}, we see that~$\tilde g$ is negative for all~\mbox{$x\geq 1$}, and thereby~$g$ is monotonously decreasing.}

\new{
In order to do so, we use the formula~\mbox{$(x\,I_1)^{\prime}(x) = x\,I_0(x)$}, see~\cite[9.6.28]{abst}, to compute the derivative
\begin{align*}
\tilde g^\prime(x)
&=
-1 - x\,I_0(x) + I_0(x) + (x-1)\,I_1(x) \\
&= 
(x-1)[I_1(x)-I_0(x)]-1 \,.
\end{align*}
Since~\mbox{$I_1(x) \leq I_0(x)$} for all~\mbox{$x \geq 0$} by~\cite[(2.3)]{Ba10} we have~\mbox{$\tilde g^\prime(x) < 0$} for all~\mbox{$x \geq 1$}, which completes the proof. \hfill \qedsymbol
\medskip
}

\begin{Theorem}
\label{Theorem:approxerrorcKB}
Assume that~\mbox{$f\in L^2(\R) \cap C(\R)$} is bandlimited with bandwidth~\mbox{$\delta \in \big(0,\,\frac{m-1}{m}\pi\big]$}. Further let~\mbox{$\varphi_{\mathrm{cKB}}$} be the continuous Kaiser--Bessel window function~\eqref{eq:varphicKB} with shape parameter~\mbox{$\beta = m\,(\pi - \delta)$} and let~\mbox{$m \in {\mathbb N} \setminus \{1\}$} be given.
\\
Then the continuous Kaiser--Bessel regularized Shannon sampling formula satisfies the error estimate
\begin{align}
\label{eq:err_cKB}
\big\| f - R_{\mathrm{cKB},m}f \big\|_{C_0(\R)} \le  \bigg(\frac{7}{8}\,m\,(\pi - \delta) + \frac{7}{\pi}\,m^2 (\pi - \delta)^2\bigg)\,\e^{-m\,(\pi-\delta)}\,\| f\|_{L^2(\R)}\,.
\end{align}
\end{Theorem}

\emph{Proof}. 
(i) 
Since~\mbox{$\varphi_{\mathrm{cKB}}$} in~\eqref{eq:varphicKB} is compactly supported on~\mbox{$[-m,\,m]$} and~\mbox{$\varphi_{\mathrm{cKB}}(m) = 0$}, we have~\mbox{$E_2(m) = 0$}. Thus, according to Theorem~\ref{Theorem:approxerror}, the approximation error can be estimated by
\begin{align*}
\big\| f - R_{\mathrm{cKB},m}f \big\|_{C_0(\R)} \le \| f\|_{L^2(\R)}\, \max_{\omega \in [-\delta,\delta]} \big| \Delta_{\mathrm{cKB}}(\omega)\big|
\end{align*}
where
\begin{align}
\label{eq:DeltacKB}
\Delta_{\mathrm{cKB}}(\omega) \coloneqq 1 - \frac{1}{\sqrt{2\pi}}\,\int_{\omega - \pi}^{\omega+\pi} {\hat \varphi}_{\mathrm{cKB}}(\tau)\,\mathrm{d}\tau\,, \quad \omega \in [-\delta,\, \delta]\,.
\end{align}
Following~\cite[p.~3, 1.1, and p.~95, 18.31]{Ob90}, the Fourier transform of~\eqref{eq:varphicKB} has the form
\begin{align}
\label{eq:hatvarphicKB}
{\hat \varphi}_{\mathrm{cKB}}(\tau) = \frac{m\, \sqrt 2}{(I_0(\beta) - 1)\,\sqrt \pi}\cdot \begin{cases} \bigg(\frac{\sinh \big(\beta \,\sqrt{1 - \nu^2}\big)}{\beta\,\sqrt{1 - \nu^2}} - \frac{\sin(\beta \nu)}{\beta \nu}\bigg) &\colon |\nu| < 1\,,\\[2ex]
\bigg(\frac{\sin \big(\beta\, \sqrt{\nu^2 - 1}\big)}{\beta\,\sqrt{\nu^2 - 1}} - \frac{\sin(\beta \nu)}{\beta \nu}\bigg) &\colon |\nu| > 1\,,
\end{cases}
\end{align}
with the scaled frequency~\mbox{$\nu = \frac{m}{\beta}\,\tau$}.
Substituting~\mbox{$\tau = \frac{\beta}{m}\,\nu$} in the integral in~\eqref{eq:DeltacKB}, the function~\mbox{$\Delta_{\mathrm{cKB}}$} reads as
\begin{align}
\label{eq:DeltacKB_subst}
\Delta_{\mathrm{cKB}}(\omega) = 1 - \frac{\beta}{m\,\sqrt{2\pi}}\,\int_{-\nu_1(-\omega)}^{\nu_1(\omega)} {\hat \varphi}_{\mathrm{cKB}}\big(\tfrac{\beta}{m}\,\nu\big)\,\mathrm{d}\nu\,, \quad \omega \in [-\delta,\, \delta]\,,
\end{align}
with the increasing linear function~\eqref{eq:nu1}.

(ii) Now we choose the shape parameter of~\eqref{eq:varphicKB} in the special form~\mbox{$\beta = m\,(\pi - \delta)$}. Thus, we have
\begin{align*}
1 = \nu_1(-\delta) \leq \nu_1(\omega) = \frac{\omega + \pi}{\pi - \delta} \leq \nu_1(\delta) = \frac{\pi + \delta}{\pi - \delta}\,, \quad \omega \in [- \delta,\, \delta]\,.
\end{align*}
In view of~\eqref{eq:hatvarphicKB} we split~\eqref{eq:DeltacKB_subst} in the form~\mbox{$\Delta_{\mathrm{cKB}}(\omega) = \Delta_{\mathrm{cKB},1} - \Delta_{\mathrm{cKB},2}(\omega)$} with
\begin{align}
\Delta_{\mathrm{cKB},1} &= 1 - \frac{2 \beta}{\pi \,(I_0(\beta) - 1)}\, \int_0^1 \bigg(\frac{\sinh\big(\beta\, \sqrt{1 - \nu^2}\big)}{\beta \,\sqrt{1 - \nu^2}} - \frac{\sin(\beta \nu)}{\beta \nu}\bigg)\,\mathrm{d}\nu\,,\notag\\
\Delta_{\mathrm{cKB},2}(\omega)  &= \frac{\beta}{\pi\,(I_0(\beta)-1)}\,\bigg(\int_1^{\nu_1(-\omega)} + \int_1^{\nu_1(\omega)}\bigg)\bigg(\frac{\sin \big(\beta\, \sqrt{\nu^2 - 1}\big)}{\beta\,\sqrt{\nu^2 - 1}} - \frac{\sin(\beta \nu)}{\beta \nu}\bigg)\, \mathrm{d}\nu\,. \label{eq:DeltacKB2}
\end{align}
Using~\cite[3.997--1]{GR80} we have
\begin{align*}
\int_0^1 \frac{\sinh\big(\beta\, \sqrt{1 - \nu^2}\big)}{\beta \,\sqrt{1 - \nu^2}}\,\mathrm{d}\nu = \frac{1}{\beta}\,\int_0^{\pi/2} \sinh(\beta\,
\cos s)\,\mathrm{d}s = \frac{\pi}{2\beta}\,{\textbf L}_0(\beta)
\end{align*}
with the modified Struve function~\eqref{eq:L0(x)}.
Additionally, by the definition of the sine integral function\new{~\eqref{eq:Si(x)}} we have
\begin{align*}
\int_0^1 \frac{\sin(\beta v)}{\beta v}\,{\mathrm d}v = \frac{1}{\beta}\,\mathrm{Si}(\beta)\,,
\end{align*}
such that we obtain
\begin{align*}
\Delta_{\mathrm{cKB},1} &= 1 - \frac{2 \beta}{\pi \,\big(I_0(\beta)-1\big)}\,\bigg(\frac{\pi}{2 \beta}\,{\textbf L}_0(\beta) - \frac{1}{\beta}\,\mathrm{Si}(\beta)\bigg) \nonumber\\
&= \frac{1}{I_0(\beta) - 1}\,\bigg(I_0(\beta) - {\textbf L}_0(\beta) - 1 + \frac{2}{\pi}\,\mathrm{Si}(\beta)\bigg)\,. %\label{eq:DeltacKB1comp}
\end{align*}
\new{Since for~\mbox{$\delta \in \big(0,\,\frac{m-1}{m}\pi\big]$} we have~\mbox{$\beta = m\,(\pi - \delta) \geq \pi$} and therefore~\mbox{$I_0(\beta)>1$}, Lemma~\ref{Lemma:abs_I0(beta)-L0(beta)}} yields
\begin{align*}
\new{|}\Delta_{\mathrm{cKB},1}\new{|} \leq \frac{1}{2\, \big(I_0(\beta) - 1 \big)}\,.
\end{align*}
Now we estimate~\mbox{$\Delta_{\mathrm{cKB},2}(\omega)$} in~\eqref{eq:DeltacKB2} for~\mbox{$\omega \in [-\delta,\, \delta]$} by the triangle inequality as
\begin{align*}
\big|\Delta_{\mathrm{cKB},2}(\omega)\big| \leq \frac{\beta}{\big(I_0(\beta)-1\big)}\,\bigg(\int_1^{\nu_1(-\omega)} + \int_1^{\nu_1(\omega)}\bigg)\bigg|\frac{\sin\big(\beta\,\sqrt{\nu^2-1}\big)}{\beta\,\sqrt{\nu^2-1}}
- \frac{\sin(\beta \nu)}{\beta \nu}\bigg|\,\mathrm{d}\nu\,.
\end{align*}
By~\cite[Lemma~4.1]{PT21a} we have
\begin{align*}
%\label{eq:|sinc-sinc|}
	\bigg|\frac{\sin\big(\beta\, \sqrt{\nu^2-1}\big)}{\beta \,\sqrt{\nu^2 - 1}} - \frac{\sin(\beta \nu)}{\beta \nu}\bigg
	| \le \frac{2}{\nu^2}\,,
	\quad \nu \ge 1\,,
\end{align*}
and therefore
\begin{align*}
	|\Delta_{\mathrm{cKB},2}(\omega)| \le \frac{4 \beta}{\pi\,\big(I_0(\beta) - 1\big)}\, \int_1^{\infty} \frac{1}{\nu^2}\,{\mathrm d}\nu = \frac{4 \beta}{\pi\,\big(I_0(\beta) - 1\big)}\,.
\end{align*}
Thereby, we conclude that
\begin{align*}
|\Delta_{\mathrm{cKB}}(\omega)| \leq \new{|}\Delta_{\mathrm{cKB},1}\new{|} + |\Delta_{\mathrm{cKB}2}(\omega)| \leq \frac{1}{I_0(\beta) -1}\,\bigg(\frac{1}{2} + \frac{4\beta}{\pi}\bigg)\,,
\quad \omega \in [-\delta,\,\delta]\,.
\end{align*}
Since by the assumption~\mbox{$0 < \delta \leq \frac{m-1}{m}\, \pi$} we have
\mbox{$\beta = m\,(\pi - \delta) \geq \pi$} for~\mbox{$m \in \mathbb N \setminus \{1\}$} and \new{by Lemma~\ref{Lemma:monotonicity} the function~\mbox{$\frac{{\mathrm e}^x}{x\,(I_0(x) - 1)}$} is monotonously decreasing for~\mbox{$x\geq 1$}}, it follows that
\begin{align*}
	\frac{{\mathrm e}^{\beta}}{\beta\,(I_0(\beta) - 1)} \leq \frac{{\mathrm e}^\pi}{\pi\,(I_0(\pi) - 1)} = 1.644967\ldots < \frac{7}{4}\,.
\end{align*}
Hence, this yields
\begin{align*}
\frac{1}{I_0(\beta) - 1}\,\bigg(\frac{1}{2} + \frac{4\beta}{\pi}\bigg)&< \frac{7\beta}{4}\,\bigg(\frac{1}{2} + \frac{4\beta}{\pi}\bigg)\,{\mathrm e}^{-\beta}
= \bigg(\frac{7}{8}\,\beta + \frac{7}{\pi}\,\beta^2\bigg)\, {\mathrm e}^{-\beta}\,.
\end{align*}
Thus, the continuous Kaiser--Bessel regularized Shannon sampling formula with the chosen shape parameter~\mbox{$\beta = m\,(\pi -\delta)$} fulfills the error estimate~\eqref{eq:err_cKB}.
This completes the proof.
\hfill\qedsymbol
\medskip

\new{Note} that already in~\cite[Theorem~4.3]{KiPoTa24} bounds on the approximation error of the Shannon sampling formula~\eqref{eq:regShannonformula} were shown for the continuous Kaiser--Bessel window function~\eqref{eq:varphicKB} with suitably chosen shape parameter~$\beta$.
Although the respective parameter~$\beta$ looks different than the one in Theorem~\ref{Theorem:approxerrorcKB}, it is basically the same, only adapted to the slightly different setting considered in~\cite{KiPoTa24}.
\new{We additionally remark that Theorem~\ref{Theorem:approxerrorcKB} finally closes the gap in our previous proof in~\cite{KiPoTa24} since Lemma~\ref{Lemma:abs_I0(beta)-L0(beta)} is a weaker version of the numerical assumption in~\cite[Figure~4.2]{KiPoTa22}, while Lemma~\ref{Lemma:monotonicity} proves the numerical assumption in~\cite[p.~23]{KiPoTa22}.}

\new{Nevertheless, similar to the previous Sections~\ref{sec:Gauss} and~\ref{sec:sinh}, the optimality of the shape parameter~\mbox{$\beta = m\,(\pi - \delta)$} for the continuous Kaiser--Bessel window function~\eqref{eq:varphicKB} is still an open problem, which could so far only be observed numerically.
\begin{Conjecture}
\label{conj:opt_cKB}
The parameter~\mbox{$\beta = m\,(\pi - \delta)$} is optimal for the Shannon sampling formula~\eqref{eq:regShannonformula} with the continuous Kaiser--Bessel window function~\eqref{eq:varphicKB}, as it guarantees the maximum decay rate of the uniform approximation error~\eqref{eq:err_approx}.
\end{Conjecture}}

\begin{Example}
\label{ex:opt_params_cKB}
Analogously as in Example~\ref{ex:opt_params_sinh} we now \new{give numerical evidence for} the optimality of the shape parameter~\mbox{$\beta = m\,(\pi - \delta)$} of the continuous Kaiser--Bessel regularized Shannon sampling formula \new{stated} in \new{Conjecture~\ref{conj:opt_cKB}}.
More precisely, for the bandlimited function~\eqref{eq:testfunction} with several bandwidth parameters~\mbox{$\delta \in \big\{\frac{\pi}{4},\, \frac{\pi}{2},\,\frac{3\pi}{4}\big\}$}, i.\,e., several oversampling rates~\mbox{$\frac{\pi}{\delta}>1$}, we consider the regularized Shannon sampling formula~\eqref{eq:regShannonformula} with the continuous Kaiser--Bessel window function~\mbox{$\varphi_{\mathrm{cKB}}$} in~\eqref{eq:varphicKB}.
The corresponding approximation error~\eqref{eq:maxerr} shall again be approximated by evaluating the given function~$f$ and its approximation~\mbox{$R_{\varphi,m}f$} at equidistant points~\mbox{$t_s\in[-1,\,1]$}, \mbox{$s=1,\dots,S$}, with~\mbox{$S=10^5$}.
To compare with the parameter \new{in Theorem~\ref{Theorem:approxerrorcKB}}, we choose the shape parameter of the continuous Kaiser--Bessel window function~\eqref{eq:varphicKB} as~\mbox{$\beta = \alpha\, m\,(\pi - \delta)$} with~\mbox{$\alpha\in\big\{\frac 12, 1, 2\big\}$}.

The outcomes for different truncation parameters~\mbox{$m\in\{2, 3, \ldots, 10\}$} are depicted in Figure~\ref{fig:error_cKB}.
As \new{expected} it can clearly be seen that the choice of~\mbox{$\alpha \neq 1$} causes worsened error decay rates with respect to~$m$.
Thus, these numerical results \new{support the proposition} that the shape parameter~\mbox{$\beta = m\,(\pi - \delta)$} of Theorem~\ref{Theorem:approxerrorcKB} is \new{indeed} optimal \new{in terms of the uniform approximation error~\eqref{eq:err_approx}}, already for very small truncation parameters~\mbox{$m\in\N\setminus\{1\}$}.
\begin{figure}[ht]
	\centering
	\captionsetup[subfigure]{justification=centering}
	\begin{subfigure}[t]{0.32\textwidth}
		\includegraphics[width=\textwidth]{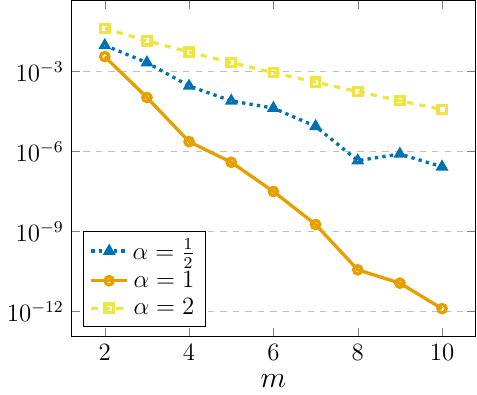}
		\caption{$\delta=\frac{\pi}{4}$}
	\end{subfigure}
	\begin{subfigure}[t]{0.32\textwidth}
		\includegraphics[width=\textwidth]{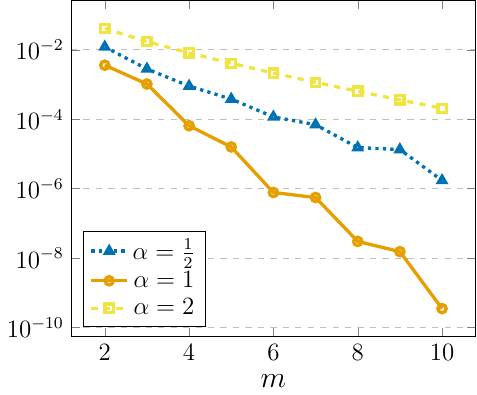}
		\caption{$\delta=\frac{\pi}{2}$}
	\end{subfigure}
	\begin{subfigure}[t]{0.32\textwidth}
		\includegraphics[width=\textwidth]{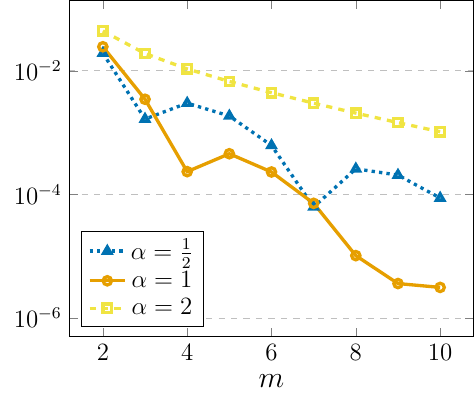}
		\caption{$\delta=\frac{3\pi}{4}$}
	\end{subfigure}
	\caption{Maximum approximation error~\eqref{eq:maxerr} using the continuous Kaiser--Bessel window function~\mbox{$\varphi_{\mathrm{cKB}}$} in~\eqref{eq:varphicKB} with different shape parameters~\mbox{$\beta=\alpha\, m (\pi-\delta)$}, \mbox{$\alpha\in\big\{\frac 12,1,2\big\}$}, for the bandlimited function~\eqref{eq:testfunction} with bandwidths~\mbox{$\delta \in \big\{\frac{\pi}{4},\frac{\pi}{2},\frac{3\pi}{4}\big\}$} and truncation parameters~\mbox{$m\in\{2, 3, \ldots, 10\}$}.
	\label{fig:error_cKB}}
\end{figure}
\end{Example}

\begin{Remark}
Note that the code files for this and all the other experiments are available on \url{https://github.com/melaniekircheis/Optimal-parameter-choice-for-regularized-Shannon-sampling-formulas}.
\end{Remark}

\section{Conclusion}

\begin{figure}[ht]
	\centering
	\captionsetup[subfigure]{justification=centering}
	\begin{subfigure}[t]{0.32\textwidth}
		\includegraphics[width=\textwidth]{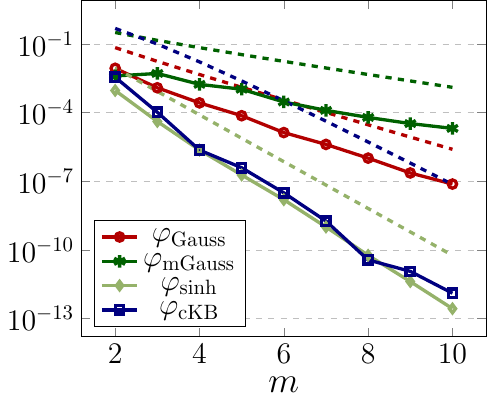}
		\caption{$\delta=\frac{\pi}{4}$}
	\end{subfigure}
	\begin{subfigure}[t]{0.32\textwidth}
		\includegraphics[width=\textwidth]{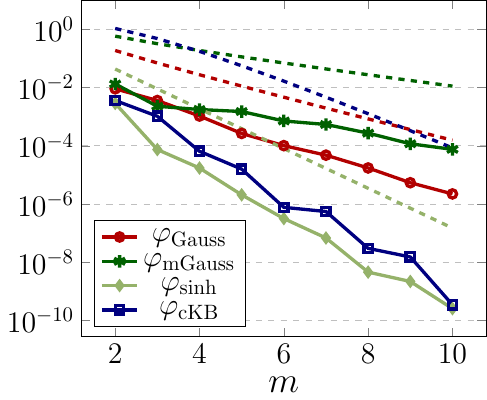}
		\caption{$\delta=\frac{\pi}{2}$}
	\end{subfigure}
	\begin{subfigure}[t]{0.32\textwidth}
		\includegraphics[width=\textwidth]{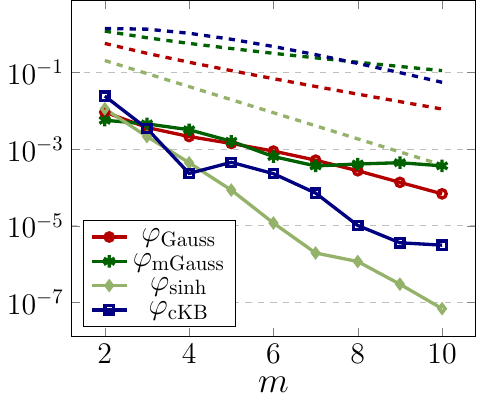}
		\caption{$\delta=\frac{3\pi}{4}$}
	\end{subfigure}
	\caption{Maximum approximation error~\eqref{eq:maxerr} (solid) and error constants (dashed) using \mbox{$\varphi\in\{\varphi_{\mathrm{Gauss}}, \varphi_{\mathrm{mGauss}}, \varphi_{\sinh}, \varphi_{\mathrm{cKB}}\}$}, see~\eqref{eq:Gaussfunction}, \eqref{eq:modGaussfunction}, \eqref{eq:varphisinh}, and~\eqref{eq:varphicKB}, for the band\-limited function~\eqref{eq:testfunction} with~\mbox{$\delta \in \big\{\frac{\pi}{4},\, \frac{\pi}{2},\,\frac{3\pi}{4}\big\}$} and~\mbox{$m\in\{2, 3, \ldots, 10\}$}.}
	\label{fig:comparison_windows}
\end{figure}

In this paper, we have studied the regularized Shannon sampling formula~\eqref{eq:regShannonformula} for the widely used Gaussian function~\eqref{eq:Gaussfunction}, the $\sinh$-type window function~\eqref{eq:varphisinh}, and the continuous Kaiser--Bessel window function~\eqref{eq:varphicKB}.
More precisely, for an arbitrary bandlimited function~\mbox{$f \in L^2(\mathbb R) \cap C(\mathbb R)$} with bandwidth~\mbox{$\delta \in (0,\, \pi)$} we have shown that the uniform approximation error~\eqref{eq:err_approx} of the regularized Shannon sampling formulas of~$f$ possess an exponential decay with respect to the truncation parameter~$m$.
In doing so, we have demonstrated that the decay rate~\mbox{$m\,(\pi-\delta)$} of the \mbox{$\sinh$-type} regularized Shannon sampling formula, see Theorem~\ref{Theorem:approxerrorsinh}, and the continuous Kaiser--Bessel regularized Shannon sampling formula, see Theorem~\ref{Theorem:approxerrorcKB}, is \new{twice as fast as} the decay rate~\mbox{$m\,(\pi-\delta)/2$} of the Gaussian regularized Shannon sampling formula, see Theorem~\ref{Theorem:approxerrorGauss}.
Note that the \mbox{$\sinh$-type} regularized Shannon sampling formula is even \new{slightly} better than the continuous Kaiser--Bessel regularized Shannon sampling formula due to the constant factors in~\eqref{eq:err_sinh} and~\eqref{eq:err_cKB}, see also Figure~\ref{fig:comparison_windows}.

\new{The main focus of this work was to elaborate on previous results in~\cite{KiPoTa22, KiPoTa24}.
First and foremost, the goal was to rigorously prove all the necessary ingredients for the proofs, thereby improving upon the results in~\cite{KiPoTa22, KiPoTa24} that lacked rigor due to the use of numerical assumptions.}
\new{In addition}, we found that the exponential decay of the approximation error of the regularized Shannon sampling formula~\eqref{eq:regShannonformula} depends highly on the shape parameter of the corresponding window function.
\new{Although the optimality} of the variance~$\sigma^2$ of the Gaussian function and of the shape parameter~$\beta$ of the \mbox{$\sinh$-type} window function and the continuous Kaiser--Bessel function is \new{still an open problem, we have strong reason to believe that our choice is best, even though this is currently only based on numerical observations.}
\new{We remark that this} further emphasizes the superiority of the \mbox{$\sinh$-type} regularized Shannon sampling formula of~$f$, since the approximation errors of the regularized Shannon sampling formulas were compared for the \new{presumably optimal} shape parameters each.

\section*{Acknowledgments}
Melanie Kircheis gratefully acknowledges the support of the BMBF grant 01$\mid$S20053A (project SA$\ell$E) and the Deutsche Forschungsgemeinschaft (DFG, German Research Foundation) – Project-ID 519323897.

Furthermore, we would like to express our sincere gratitude to the reviewers and the editor for their careful reading, valuable comments, and substantial suggestions, all of which helped us improve this work considerably.

\end{document}